\documentclass{article}
\usepackage{amssymb}
\usepackage{amsfonts}
\usepackage{graphicx}
\usepackage{amsmath}
\usepackage{booktabs}
\usepackage{multirow}

\newtheorem{algorithm}{Algorithm}

\newtheorem{lemma}{Lemma}

\newtheorem{proposition}{Proposition}

\newenvironment{proof}[1][Proof]{\textbf{#1.} }{\ \rule{0.5em}{0.5em}}

\begin{document}

\date{}
\title{A fast and accurate numerical approach for electromagnetic inversion}
\author{Eleonora Denich\thanks{%
Dipartimento di Matematica e Geoscienze, Universit\`{a} di Trieste,
eleonora.denich@phd.units.it} \and Paolo Novati \thanks{%
Dipartimento di Matematica e Geoscienze, Universit\`{a} di Trieste,
novati@units.it} \and Stefano Picotti \thanks{Department of Geophisics, Istituto Nazionale di Oceanografia e di Geofisica Sperimentale (OGS), Borgo Grotta Gigante 42c, 34010 Sgonico, Trieste, Italy, spicotti@inogs.it}}
\maketitle

\begin{abstract}
This paper deals with the solution of Maxwell's equations to model the electromagnetic fields in the case of a layered earth. The integrals involved in the solution are approximated by means of a novel approach based on the splitting of the reflection term. The inverse problem, consisting in the computation of the unknown underground conductivity distribution from a set of modeled magnetic field components, is also considered. Two optimization algorithms are applied, based on line- and global-search methods, and a new minimization approach is presented.  Several EM surveys from the ground surface are simulated, considering the horizontal coplanar (HCP) and perpendicular (PRP) magnetic dipolar configurations. The numerical experiments, carried out for the study of river-levees integrity, allowed to estimate the errors associated to these kind of investigations, and confirm the reliability of the technique.
\end{abstract}

%\begin{keyword}
%Maxwell equations \sep Gauss-Kronrod quadrature \sep EM response \sep EM inversion
%\MSC[2010] 78A25 \sep  78A46 \sep 65D32 \sep 78M50
%\end{keyword}

\section{Introduction}

The aim of electromagnetic (EM) sounding methods in geophysics is to obtain information about the subsoil conductivity distribution, by recorded measurements of the EM field components. One technique consists in placing magnetic dipoles above the surface, composed of a transmitter coil and different couples of adjacent receiver coils. The receiver couples (i.e. dual coils) are placed at different distances (offsets) from the transmitter coil. The most common dipole geometry consists of transmitter and receiver loops that can be horizontal co-planar (HCP configuration) and perpendicular (PRP configuration). 
Among the available instruments, the DUALEM (DUAL-geometry Electro-Magnetic; http://www.dualem.com) system is often used. The receiver couples are placed at $2,4,6$ and $8$ m from the transmitter coil, and typical source-receiver geometries are the horizontal coplanar (HCP configuration) and perpendicular (PRP configuration) coils.
The electromagnetic induction effect, encoded in the first-order linear differential Maxwell equations, produces eddy
alternating currents in the soil which on their turn, induce response
EM fields, that can be used to determine the conductivity profile of the ground.

For simplicity, the local subsurface structures below the DUALEM instrument can be assumed composed by horizontal and homogeneous layers. 
Under this assumption, general integral solutions of Maxwell equations (i.e., the EM fields) for vertical and horizontal magnetic dipoles, can be derived (see \cite{WH}) and represented as Hankel transforms of order $l$, as follows: 
\begin{equation}
F(r,p)=\int_{0}^{\infty }f(\lambda ,p)J_{l}(\lambda r)d\lambda ,
\label{Hankel}
\end{equation}%
where $r$ is the offset and $J_{l}$ is the Bessel function of the first kind of order $l$. The function $f(\lambda,p)$ contains the vector $p$ characterizing the subsurface model parameters, i.e. the conductivity and the thickness of each layer. As for the Bessel functions we refer here to \cite{Bessel} for an overview.
In the case of conductivities of geological materials, only the imaginary part of the complex function $F(r,p)$ is considered. Indeed, while the instruments (e.g., the DUALEM) used for experimental surveys collect both in-phase and quadrature measurements of the magnetic field components, the former become important only over highly conductive materials, and are particularly effective for locating confined conductors, such as metal bodies, or boulders of graphite or sulfide.

By a mathematical point of view, in general, the function $F(r,p)$ is not analitically computable and therefore it is necessary to employ a numerical scheme. Anyway, the slowly decay of the oscillations determined by the Bessel function makes the problem very difficult to handle, because traditional quadrature rules typically fail to converge. In order to numerically evaluate these kind of integrals, in 1971 Ghosh \cite{Ghosh} introduced the digital filtering algorithm for the direct interpretation of geoelectrical resistivity soundings measuraments. This method is essentially a standard quadrature rule, but the main difference is that the weights are computed by solving a linear equation obtained by imposing the rule to be correct on a set of training functions (not polynomials) for which the corresponding integral (\ref{Hankel}) is known (e.g., \cite{G-S}).
%For example, in \cite{G-S} the following Hankel transforms are considered: 
%\begin{align*}
%& \int_{0}^{\infty }\lambda e^{-c\lambda}J_{0}(\lambda r)d\lambda =\frac{c%
%}{(c^{2}+r^{2})^{3/2}}, \\
%& \int_{0}^{\infty }(\lambda e^{-c\lambda }+\alpha \lambda ^{2}e^{-c\lambda
%^{2}})J_{1}(\lambda r)d\lambda =\frac{r}{(c^{2}+r^{2})^{3/2}+\alpha \frac{%
%re^{-(r^{2}/4c)}}{4c^{2}}},
%\end{align*}%
%with $\alpha ,r,c>0$. The set of training functions is obtained by moving $\alpha$ and $c$.

In general, considering the change of variables $r=e^{x}$ and $\lambda =e^{-y}$, equation (\ref{Hankel}) becomes
\begin{equation*}
e^x F(e^x,p)= \int_{-\infty}^{+\infty} f(e^{-y},p) J_l(e^{x-y}) dy,
\end{equation*}
and can be rewritten as the convolution
\begin{equation*}
\tilde{F}(x,p) =\int_{-\infty}^{+\infty} \tilde{f}(y,p) \tilde{J}_l(x-y) dy = \int_{-\infty}^{+\infty} \tilde{f}(x-y,p) \tilde{J}_(y) dy,
\end{equation*}
where $\tilde{F}(x,p)=e^xF(e^x,p)$, $\tilde{f}(y,p)=f(e^{-y},p)$ and $\tilde{J}_l(x-y)=J_l(e^{x-y})$. Then, given a set of training functions $\left\{ f_{k}\right\} _{k=-N,...,N}$, the method
consists in setting $2N+1$ grid points $\lambda _{j}=e^{-jy}$, $j=-N,...,N$, and then
prescribes to solve with respect to the filter coefficients $w_{j}$, $%
j=-N,...,N$, the linear system%
\begin{equation}
\tilde{F}_{k}(x,p)=\sum\nolimits_{j=-N}^{N}\tilde{f}_{k}(x-y,p)w_{j},\quad
k=-N,...,N.  \label{convolution_equation}
\end{equation}%
Finally,
the resulting $2N+1$-points digital linear filter approximation of equation (%
\ref{Hankel}) is 
\begin{equation}
F(r,p)\approx \frac{1}{r}\sum\nolimits_{j=-N}^{N}f\left( \frac{\lambda _{j}}{%
r},p\right) w_{j}.  \label{df}
\end{equation}

Later, this algorithm has been reconsidered by many authors (see e.g. \cite{anderson2,JS,KGP,anderson1,G-S,Emdpler,Guptasarma,Werthmuller}) to handle specific EM problems involving integrals of type (\ref{Hankel}), and improvements to the determination of filter coefficients have been developed (see \cite{WienerHopf,Kong}). In particular, Koefoed and Dirks \cite{WienerHopf}
proposed a Wiener-Hopf least-squares method, while Kong \cite{Kong} uses the GMRES, to solve the linear system (\ref{convolution_equation}).

In this work we consider a different approach that is much simpler and reliable for integrals of type (\ref{Hankel}). Indeed, in this framework the function $f(\lambda, p )$ can be decomposed as $f_{1}(\lambda, p )+f_{2}(\lambda, p )$ in which $f_{1}(\lambda,p )$ is such that $F_1(r, p)$ is known exactly, and the oscillating function $f_{2}(\lambda, p )$ decays exponentially. For realistic sets of parameters, the oscillations are quite rapidly damped and the corresponding integral $F_2(r, p)$
can be accurately computed by a classical quadrature rule on finite
intervals. We remark that the frequency of the oscillations increases with $r$, so that the integration of $F(r,p)$ becomes more difficult. In this view the splitting allows to reduce the sensitivity with respect to $r$. The idea of splitting the integral in two parts has originally
been introduced in \cite{CR1} where, however, the second integral is still
approximated as in (\ref{df}) by digital filtering, without exploiting the
fast decay of the oscillations.

In order to reconstruct realistic subsurface structures from EM measurements, a tomographic approach is needed.
A typical inversion approach consists in an iterative procedure involving the computation of the EM response of a layered model (forward modelling) and the solution of the inverse problem. The algorithm attempts to minimize the mismatch between the measured data and the predicted data, by updating the model parameters at each iteration.
Having at disposal a reliable method for evaluating (\ref{Hankel}), here we
also consider the inverse problem of computing the model parameters (i.e., conductivity and thickness of the layers) from a set of measured field values at different offsets. To this
purpose we employ two optimization algorithms based on the BFGS line-search method (\cite{Broyden,Fletcher,Goldfarb,Shanno}) and on the Simulated Annealing (SA) global-search technique (\cite{Kirkpatrick,Goffe}). In order to reduce as much as possible the number of integral evaluations, we also derive an analytic approximation of these integrals that can be used 
in the initial iterations of the tomographic procedure, to have a first estimate of the solution.
% as a preprocessing phase to locate with a good precision the set of parameters.

This paper is organized as follows. In Section 2 we define the integral formulations of the EM fields. Section 3 deals with the numerical approach used for the computation of the integrals. In Section 4 we derive useful approximations in the case of a layered earth. In Section 5 we introduce the inverse problem and propose an optimization algorithm to reconstruct the underground conductivity distribution. Finally, in Section 6 we present some realistic examples of river levees conductivity models, and report the numerical results in the case of a $3$-layered subsoil.

\section{Definition of the forward model}

\begin{figure}[t]
\begin{center}
\includegraphics[scale=0.6]{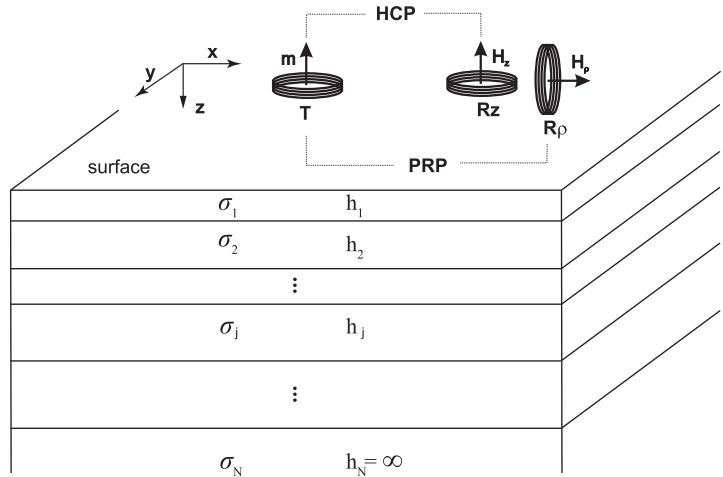}
\end{center}
\caption{Layered underground conductivity model and dual-coil configuration inside the DUALEM instrument.}
\label{strati}
\end{figure}

Consider a layered underground model as in Figure \ref{strati}, where $%
\sigma _{j}$ and $h_{j}$, $j=1,...,N$, represent conductivity and thickness
of the $j$-th layer, respectively. The deeper layer is assumed to have infinite
thickness. Let $f$ and $m$ be the transmitter's frequency and
the magnetic moment, $\mu $ the magnetic permeability of vacuum and $r$ the offset. The angular frequency is $\omega=2 \pi f$. In order to derive the mathematical formulation of the fields, it is convenient to introduce a
cylindrical polar coordinate system $(\rho ,\phi ,z)$, with the longitudinal
axis downward directed, such that the ground plane coincides with the plane $%
z=0$, the transmitter is located at the origin and the receiver is placed
along the polar axis.

Let us consider the DUALEM instrument, for which the dipole geometry consists of
a transmitter loop (T) and many dual receiver loops (Rz, R$\rho$) that are horizontal co-planar (HCP configuration)
and perpendicular (PRP configuration). Figure \ref{strati} represents an example of dual-coil configuration inside the DUALEM, where the offset is 2 m and 2.1 m between T and the first couple of receivers Rz, R$\rho$, respectively. Other couples of receivers are located at 4, 6 and 8 m offset. 
In this setting, the theoretical components of the magnetic field at the receiver location on the surface are given by (see \cite{WH} and \cite{CR1}) 
\begin{align}
H_{z}^{(N)}& =\frac{m}{4\pi }\int_{0}^{\infty }(1+R_{0}(\lambda ))\lambda
^{2}J_{0}(\lambda r)d\lambda ,  \label{Hz} \\
H_{\rho }^{(N)}& =\frac{m}{4\pi }\int_{0}^{\infty }(1-R_{0}(\lambda
))\lambda ^{2}J_{1}(\lambda r)d\lambda ,  \label{Hr}
\end{align}%
for HCP and PRP configurations, respectively. In the above formulas $%
R_{0}(\lambda )$ is the reflection term, recursively defined by 
\begin{equation} \label{x}
\begin{split}
R_{0}(\lambda )& =\frac{R_{1}(\lambda )+\Psi _{1}(\lambda )}{R_{1}(\lambda
)\Psi _{1}(\lambda )+1}, \\
R_{j}(\lambda )& =\frac{R_{j+1}(\lambda )+\Psi _{j+1}(\lambda )}{%
R_{j+1}(\lambda )\Psi _{j+1}(\lambda )+1}e^{-2u_{j}(\lambda )h_{j}},\quad
j=1,...,N-1, \\
R_{N}(\lambda )& =0,
\end{split}
\end{equation}%
where%
\begin{equation*}
\Psi _{j}(\lambda )=\frac{u_{j-1}(\lambda )-u_{j}(\lambda )}{u_{j-1}(\lambda
)+u_{j}(\lambda )},\quad j=1,...,N,
\end{equation*}%
in which $u_{0}(\lambda )=\lambda $ and $u_{j}(\lambda )=\sqrt{\lambda
^{2}-k_{j}^{2}}$, $k_{j}=\sqrt{-i\omega \mu \sigma _{j}}$, for $j=1,\ldots ,N
$. In order to avoid redundancies and to reduce the length of some formulas,
in the sequel we simply write $R_{j}$, $\Psi _{j}$, $u_{j}$, in place of $%
R_{j}(\lambda )$, $\Psi _{j}(\lambda )$, $u_{j}(\lambda )$.

We observe that if $N=1$, representing the situation of homogeneous earth, then $R_{1}=0$, $R_{0}=\Psi _{1}$, and the
integrals (\ref{Hz}) e (\ref{Hr}) can be exactly evaluated. Indeed, in \cite%
{WH} it is shown that in this situation 
\begin{align*}
H_{z}^{(1)}& =\frac{m}{2\pi r}\frac{\partial }{\partial r}\left(
r\int_{0}^{\infty }\frac{\lambda ^{2}}{\lambda +u_{1}}J_{1}(\lambda
r)d\lambda \right)  \\
& =\frac{m}{2\pi r}\frac{\partial }{\partial r}\left( -r\frac{\partial }{%
\partial r}\int_{0}^{\infty }\frac{\lambda }{\lambda +u_{1}}J_{0}(\lambda
r)d\lambda \right)  \\
& =\frac{m}{2\pi r}\frac{\partial }{\partial r}\left( -r\frac{\partial }{%
\partial r}\int_{0}^{\infty }\frac{\lambda (\lambda -u_{1})}{\lambda
^{2}-u_{1}^{2}}J_{0}(\lambda r)d\lambda \right)  \\
& =\frac{m}{2\pi rk_{1}^{2}}\frac{\partial }{\partial r}\left( -r\frac{%
\partial }{\partial r}\left( \int_{0}^{\infty }\lambda ^{2}J_{0}(\lambda
r)d\lambda -\int_{0}^{\infty }\lambda u_{1}J_{0}(\lambda r)d\lambda \right)
\right) , \\
H_{\rho }^{(1)}& =\frac{m}{4\pi }\frac{\partial }{\partial r}%
\int_{0}^{\infty }\frac{\lambda -u_{1}}{\lambda +u_{1}}\lambda J_{0}(\lambda
r)d\lambda  \\
& =\frac{m}{4\pi }\frac{\partial }{\partial r}\int_{0}^{\infty }\frac{%
\lambda -u_{1}}{\lambda +u_{1}}\frac{\lambda +u_{1}}{\lambda +u_{1}}\lambda
J_{0}(\lambda r)d\lambda  \\
& =\frac{mk_{1}^{2}}{4\pi }\frac{\partial }{\partial r}\int_{0}^{\infty }%
\frac{\lambda }{(\lambda +u_{1})^{2}}J_{0}(\lambda r)d\lambda .
\end{align*}%
An analytical expression for the above integrals can be achieved by using \cite[p.707]%
{GR} and \cite[Vol.2, pp.8-9]{Erdelyi} respectively, to obtain 
\begin{align*}
H_{z}^{(1)}& =\frac{m}{2\pi k_{1}^{2}r^{5}}\left[
9-(9+9ik_{1}r-4k_{1}^{2}r^{2}-ik_{1}^{3}r^{3})e^{-ik_{1}r}\right] , \\
H_{\rho }^{(1)}& =-\frac{mk_{1}^{2}}{4\pi r}\left[ I_{1}\left( \frac{ik_{1}r%
}{2}\right) K_{1}\left( \frac{ik_{1}r}{2}\right) -I_{2}\left( \frac{ik_{1}r}{%
2}\right) K_{2}\left( \frac{ik_{1}r}{2}\right) \right] ,
\end{align*}%
where $I_{l}$ e $K_{l}$ are the modified Bessel functions of order $l$, of
the first and second type respectively.

Defining, for $l=0,1$
\begin{align}
g_{l}\left( \lambda \right) &=\left( R_{0}-\Psi _{1}\right) \lambda
^{2}J_{l}(\lambda r), \label{g_l} \\
q_l(\lambda)&= [1+(-1)^l \Psi_1] \lambda^2 J_l(\lambda r),
\end{align}%
we can write (see (\ref{Hz}), (\ref{Hr}))
\begin{align*}
(1+R_0)\lambda^2 J_0(\lambda r) &= g_0(\lambda) + q_0(\lambda), \\
(1-R_0)\lambda^2 J_1(\lambda r) &= q_1(\lambda)- g_1(\lambda).
\end{align*}
Since $q_l(\lambda)$ represents the integrand function in the case $N=1$ ($R_1=0$ and hence $g_l(\lambda)=0$), we have that
\begin{align*}
\frac{m}{4 \pi} \int_0^{\infty} q_0(\lambda) d \lambda &= H_z^{(1)}, \\
\frac{m}{4 \pi} \int_0^{\infty} q_1(\lambda) d \lambda &= H_{\rho}^{(1)},
\end{align*}
and therefore
\begin{align}
H_{z}^{(N)} & =\frac{m}{4\pi }\int_{0}^{\infty }g_{0}\left( \lambda \right) d\lambda
+H_{z}^{(1)}, \label{eq9} \\
H_{\rho}^{(N)} & =-\frac{m}{4\pi }\int_{0}^{\infty }g_{1}\left( \lambda \right) d\lambda
+H_{\rho}^{(1)}. \label{eq10}
\end{align}

\section{Numerical computation of the fields}

Our strategy for the computation of $H_{z}^{(N)}$ and $H_{\rho }^{(N)}$ in the general case of $N$ layers is based on the observation that the functions $g_l(\lambda)$, $l=0,1$, exponentially decay with respect to $\lambda$. To prove this behavior, some preliminary results are necessary.

\begin{lemma} \label{lemmaA}
It holds 
\begin{equation*}
R_0-\Psi_1 = R_1 \left( 1+ \mathcal{O} \left( \frac{1}{\lambda^2} \right) \right), \quad {\rm as} \; \lambda \rightarrow + \infty.
\end{equation*}
\end{lemma}
\begin{proof}
Using the definitions of $R_0$ and $\Psi_1$ we have
\begin{equation}
R_{0}-\Psi _{1}=\frac{4R_{1}u_{1}}{R_{1}k_{1}^{2}+(\lambda +u_{1})^{2}}%
\lambda.   \label{diff}
\end{equation}
Since \begin{equation*}
u_1= \lambda \sqrt{1-\frac{k_1^2}{\lambda^2}}= \lambda \left( 1-\frac{k_1^2}{2\lambda^2}+ \mathcal{O} \left(\frac{1}{\lambda^4} \right) \right),
\end{equation*}
as $\lambda \rightarrow +\infty$, we can write
\begin{equation*}
\begin{split}
R_0-\Psi_1 &=\frac{4 R_1 u_1 \lambda}{R_1 k_1^2 + (\lambda + u_1)^2} = \frac{4 R_1 \lambda^2 \left( 1-\frac{k_1^2}{2\lambda^2}+\mathcal{O} \left( \frac{1}{\lambda^4} \right) \right)}{R_1 k_1^2 + \lambda^2 \left( 2-\frac{k_1^2}{2\lambda^2}+\mathcal{O} \left( \frac{1}{\lambda^4} \right) \right)^2}\\
&= R_1 \left( 1+ \mathcal{O} \left( \frac{1}{\lambda^2} \right) \right).
\end{split}
\end{equation*}
\end{proof}

\begin{lemma}
The function $R_1$ can be written in the following form
\begin{equation}
R_{1}=\sum\nolimits_{k=1}^{N-1}\delta _{k}(\lambda )\exp \left(
-\gamma _{k}\lambda \right),  \label{dexp}
\end{equation}%
where
\begin{align}
\gamma_{k}&=\sum\nolimits_{i=1}^{k}c_{i}> 0, \quad k=1, \ldots, N-1, \notag \\
c_{i}&=2h_{i}\sqrt{1+\frac{i\omega \mu \sigma _{i}}{\lambda ^{2}}},\quad
i=1,\ldots ,N-1. \label{E}
\end{align}
\end{lemma}
\begin{proof}
In order to demonstrate (\ref{dexp}), starting from $j=1$, by induction we show that 
\begin{equation} \label{induzione}
R_{N-j} = \sum_{k=1}^j \delta_k^{(N-j)} \exp(-(\gamma_{N+k-(j+1)}-\gamma_{N-(j+1)}) \lambda), \quad {\rm for} \; j=1, \ldots, N-1,
\end{equation}
in which $\gamma_0=0$.
We first observe that
\begin{equation*}
2u_{i}h_{i}=2h_{i}\sqrt{\lambda ^{2}+i\omega \mu \sigma _{i}}=2h_{i}\sqrt{1+ \frac{i\omega \mu \sigma _{i}}{\lambda ^{2}}}\cdot \lambda =c_{i}\lambda.
\end{equation*}
Let $j=1$. By (\ref{x}) and defining 
\begin{equation} \label{G}
\delta_{1}^{(N-1)}=\Psi _{N}=\frac{\sqrt{\lambda^2-k_{N-1}^2}-\sqrt{\lambda^2-k_{N}^2}}{\sqrt{\lambda^2-k_{N-1}^2}+\sqrt{\lambda^2-k_{N}^2}},
\end{equation}
we have that (\ref{induzione}) is correct for $j=1$ because $R_N=0$. Assuming that (\ref{induzione}) is also correct for a given $j<N-1$, by (\ref{x})
\begin{equation} \label{passo_induzione}
\begin{split}
R_{N-(j+1)}& =\frac{R_{N-j}+\Psi _{N-j}}{R_{N-j}\Psi _{N-j}+1}%
e^{-c_{N-(j+1)}\lambda } \\
%& =\delta_1^{(N-(j+1))}e^{-c_{N-(j+1)}\lambda}+\delta_2^{(N-(j+1))}e^{-(c_{N-j}+c_{N-(j+1)})\lambda} +\ldots \\
%&+\delta_j^{(N-(j+1))}e^{-(c_{N-2}+\ldots+c_{N-(j+1)})\lambda}+\delta_{j+1}^{(N-(j+1))}e^{-(c_{N-1}+\ldots + c_{N-(j+1)})\lambda}\\
&=\sum_{k=1}^{j+1} \delta_k^{(N-(j+1))} \exp(-(\gamma_{N+k-(j+2)}-\gamma_{N-(j+2)}) \lambda),
\end{split}
\end{equation}
in which we have defined 
\begin{equation} \label{T}
\begin{split}
\delta_1^{(N-(j+1))}&=\frac{\Psi_{N-j}}{R_{N-j}\Psi_{N-j}+1},\\
\delta_k^{(N-(j+1))}&=\frac{\delta_{k-1}^{(N-j)}}{R_{N-j}\Psi_{N-j}+1}, \quad k=2, \ldots, j+1.
\end{split}
\end{equation}
% the term $R_{1}$ can be written using the recursion
%\begin{align*}
%R_{N-1}& =\delta _{1}^{(N-1)}e^{-c_{N-1}\lambda }, \\
%R_{N-2}& =\frac{R_{N-1}+\Psi _{N-1}}{R_{N-1}\Psi _{N-1}+1}e^{-c_{N-2}\lambda
%}=\frac{\delta _{1}^{(N-1)}e^{-c_{N-1}\lambda }+\Psi_{N-1}}{\Psi _{N-1}\delta
%_{1}^{(N-1)}e^{-c_{N-1}\lambda }+1}e^{-c_{N-2}\lambda } \\
%& =\frac{\delta _{1}^{(N-1)}}{\Psi _{N-1}\delta
%_{1}^{(N-1)}e^{-c_{N-1}\lambda }+1}e^{-(c_{N-1}+c_{N-2})\lambda }+\frac{\Psi
%_{N-1}}{\Psi _{N-1}\delta _{1}^{(N-1)}e^{-c_{N-1}\lambda }+1}%
%e^{-c_{N-2}\lambda } \\
%& =\delta _{2}^{(N-2)}e^{-(c_{N-1}+c_{N-2})}\lambda +\delta
%_{1}^{(N-2)}e^{-c_{N-2}\lambda }, \\
%\vdots &  \\
%R_{N-k}& =\frac{R_{N-(k-1)}+\Psi _{N-(k-1)}}{R_{N-(k-1)}\Psi _{N-(k-1)}+1}%
%e^{-c_{N-k}\lambda } \\
%& =\delta _{k}^{(N-k)}e^{-(c_{N-1}+\ldots +c_{N-k})\lambda }+\delta
%_{k-1}^{(N-k)}e^{-(c_{N-2}+\ldots +c_{N-k})\lambda }+\ldots  \\
%& +\delta _{2}^{(N-k)}e^{-(c_{N-(k-1)}+c_{N-k})\lambda }+\delta
%_{1}^{(N-k)}e^{-c_{N-k}\lambda }, \\
%\vdots &  \\
%R_{1}& =\delta _{N-1}^{(1)}e^{-(c_{N-1}+\ldots +c_{1})\lambda }+\delta
%_{N-2}^{(1)}e^{-(c_{N-2}+\ldots +c_{1})\lambda }+\ldots +\delta
%_{2}^{(1)}e^{-(c_{2}+c_{1})\lambda }+\delta _{1}^{(1)}e^{-c_{1}\lambda }.
%\end{align*}%
Finally, setting $\delta_j=\delta_j^{(1)}$, $j=1, \ldots, N-1$, we obtain the result.
%\begin{align*}
%R_{1}& =\delta _{1}e^{-\gamma \lambda }+\delta
%_{2}e^{-\gamma _{2}\lambda }+\ldots +\delta _{N-2}e^{-\gamma _{N-2}\lambda
%}+\delta _{N-1}e^{-\gamma _{N-1}\lambda } \\
%& =\sum\nolimits_{j=1}^{N-1}\delta _{j}(\lambda )\exp \left( -\gamma
%_{j}\lambda \right) 
%\end{align*}
\end{proof}

\begin{lemma} \label{lemma}
For $\lambda \rightarrow + \infty$,
\begin{equation}
R_1=\sum_{k=1}^{N-1} \frac{k_{k+1}^2-k_{k}^2}{4 \lambda^2} \left( 1+ \mathcal{O} \left(\frac{1}{\lambda^2} \right) \right) \exp(-\gamma_{k} \lambda).
\end{equation}
\end{lemma}
\begin{proof}
Starting from $j=1$, by induction we show that
\begin{equation} \label{ind_asterisco}
R_{N-j}=\sum_{k=1}^{j} \frac{k_{N+k-j}^2-k_{N+k-(j+1)}^2}{4 \lambda^2} \left( 1+ \mathcal{O} \left(\frac{1}{\lambda^2} \right) \right) \exp(-(\gamma_{N+k-(j+1)}-\gamma_{N-(j+1)}) \lambda), 
\end{equation}
for $j=1,\ldots, N-1$. We observe that, for $i=1, \ldots, N$
\begin{align*}
\Psi_i &= \frac{\sqrt{\lambda^2-k_{i-1}^2}-\sqrt{\lambda^2-k_{i}^2}}{\sqrt{\lambda^2-k_{i-1}^2}+\sqrt{\lambda^2-k_{i}^2}} \\
&= \frac{k_i^2-k_{i-1}^2}{\left[ \lambda \left( 1+\frac{k_{i-1}^2}{2 \lambda^2}+ \mathcal{O} \left( \frac{1}{\lambda^4} \right) \right) + \lambda \left( 1+\frac{k_{i-1}^2}{2 \lambda^2}+ \mathcal{O} \left( \frac{1}{\lambda^4} \right) \right) \right]^2} \\
&= \frac{k_i^2-k_{i-1}^2}{4 \lambda^2 \left( 1+ \mathcal{O} \left( \frac{1}{\lambda^2} \right) \right)} = \frac{k_i^2-k_{i-1}^2}{4 \lambda^2} \left( 1+ \mathcal{O} \left( \frac{1}{\lambda^2} \right) \right).
\end{align*}
By (\ref{induzione}) with $j=1$ and (\ref{G}), we have that
\begin{equation*}
\begin{split}
R_{N-1} &= \delta_1^{(N-1)} \exp(-(\gamma_{N-1}-\gamma_{N-2})\lambda) \\
&= \frac{k_N^2-k_{N-1}^2}{4\lambda^2} \left(1+\mathcal{O}\left(\frac{1}{\lambda^2} \right) \right) \exp(-(\gamma_{N-1}-\gamma_{N-2})\lambda). 
\end{split}
\end{equation*}
%where the last equality comes from 
%\begin{equation} \label{BIG}
%\left( 1+ \mathcal{O} \left( \frac{1}{\lambda^2} \right) \right)^{-\tilde{\gamma}_k \lambda} = 1+ \mathcal{O} \left( \frac{1}{\lambda} \right).
%\end{equation}
Therefore, (\ref{ind_asterisco}) holds true for $j=1$. Assuming that (\ref{ind_asterisco}) is also correct for a given $j<N-1$, by (\ref{passo_induzione}) we have
\begin{equation*}
R_{N-(j+1)}= \sum_{k=1}^{j+1} \delta_{k}^{(N-(j+1))} \exp(-(\gamma_{N+k-(j+2)}-\gamma_{N-(j+2)}) \lambda),
\end{equation*}
where by (\ref{T}) and using the induction hypothesis on $R_{N-j}$
\begin{equation*}
\begin{split}
\delta_1^{(N-(j+1))}&=\frac{\Psi_{N-j}}{R_{N-j}\Psi_{N-j}+1}\\
&=\frac{\frac{k_{N-j}^2-k_{N-j-1}^2}{4 \lambda^2} \left( 1+ \mathcal{O} \left( \frac{1}{\lambda^2} \right) \right)}{\mathcal{O} \left( \frac{1}{\lambda^2} \right) \frac{k_{N-j}^2-k_{N-j-1}^2}{4 \lambda^2} \left( 1+ \mathcal{O} \left( \frac{1}{\lambda^2} \right) \right) +1} \\
&=\frac{k_{N-j}^2-k_{N-j-1}^2}{4 \lambda^2} \left(1+\mathcal{O}\left(\frac{1}{\lambda^2} \right) \right), \\
\delta_{k}^{(N-(j+1))}&=\frac{\delta_{k-1}^{(N-j)}}{R_{N-j}\Psi_{N-j}+1} \\
&= \frac{\frac{k_{N+k-(j+1)}^2-k_{N+k-(j+2)}^2}{4 \lambda^2} \left(1+\mathcal{O}\left(\frac{1}{\lambda^2} \right) \right)}{\mathcal{O} \left( \frac{1}{\lambda^2} \right) \frac{k_{N+k-(j+1)}^2-k_{N+k-(j+2)}^2}{4 \lambda^2} \left(1+\mathcal{O}\left(\frac{1}{\lambda^2} \right) \right) +1} \\
&=\frac{k_{N+k-(j+1)}^2-k_{N+k-(j+2)}^2}{4 \lambda^2} \left(1+\mathcal{O}\left(\frac{1}{\lambda^2} \right) \right), \quad k=2, \ldots, j+1.
\end{split}
\end{equation*}
%\begin{equation*}
%R_{N-(j+1)} = \sum_{k=1}^{j+1} \frac{k_{N+k-(j+1)}^2-k_{N+k-(j+2)}^2}{4\lambda^2} \left(1+\mathcal{O}\left(\frac{1}{\lambda} \right) \right) \exp(-(\tilde{\gamma}_{N+k-(j+2)}-\tilde{\gamma}_{N-(j+2)}) \lambda),
%\end{equation*}
%and the Lemma holds true.
\end{proof}

Using the above Lemmas we can prove the following result, which finally expresses the asymptotic bahavior of $R_1$, and hence of $R_0 - \Psi_1$ by Lemma \ref{lemmaA}.

\begin{proposition} \label{proposizione1}
For $\lambda \rightarrow + \infty$,
\begin{equation}
R_1=\sum_{k=1}^{N-1} \frac{k_{k+1}^2-k_{k}^2}{4 \lambda^2} \exp(-\tilde{\gamma}_{k} \lambda) \left( 1+ \mathcal{O} \left(\frac{1}{\lambda} \right) \right),
\end{equation}
where $\tilde{\gamma}_k=2 \sum_{i=1}^k h_i$.
\end{proposition}
\begin{proof}
First of all, by (\ref{E}) we have
\begin{equation} \label{GAMMA}
\gamma_k = \sum_{i=1}^k c_i = 2 \sum_{i=1}^k h_i \left(1+\mathcal{O}\left(\frac{1}{\lambda^2} \right) \right)=\tilde{\gamma}_k \left(1+\mathcal{O}\left(\frac{1}{\lambda^2} \right) \right), \quad {\rm for} \; k=1, \ldots, N-1.
\end{equation}
Therefore
\begin{equation*}
\begin{split}
\exp(-\gamma_k \lambda)&= \exp \left( -\tilde{\gamma}_k \lambda \left( 1+ \mathcal{O} \left( \frac{1}{\lambda^2} \right) \right) \right) \\
&= \exp(-\tilde{\gamma}_k \lambda) \left( 1+ \mathcal{O} \left( \frac{1}{\lambda} \right) \right),
\end{split} 
\end{equation*}
where the last equality comes from 
\begin{equation*}
\exp \left( -\tilde{\gamma}_k \mathcal{O} \left( \frac{1}{\lambda} \right) \right) = 1+ \mathcal{O} \left( \frac{1}{\lambda} \right).
\end{equation*}
The result then follows straighfully from Lemma \ref{lemma}. 
\end{proof}

By Lemma \ref{lemmaA} and Proposition \ref{proposizione1}, we have that the functions $g_l$, which are the integrand functions in (\ref{eq9}) and (\ref{eq10}), can be written as
\begin{equation}
\begin{split}
g_l &= (R_0-\Psi_1) \lambda^2 J_l(\lambda r) \\
&= \sum_{k=1}^{N-1} \frac{k_{k+1}^2-k_{k}^2}{4} \exp(-\tilde{\gamma}_k \lambda) J_l(\lambda r) \left( 1+ \mathcal{O} \left( \frac{1}{\lambda} \right) \right), \label{pallino}
\end{split}
\end{equation}
so that the oscillations due to the Bessel functions are rapidely damped.

%\begin{figure}[tbp]
%\caption{The imaginary part of the function $(1- \Psi_1) \lambda^2 J_1(\lambda r)$, with $r=2$ $m$, $f= 10$ $kHz$, $m=1$ $A/m^2$
%and $\protect\sigma_1=33$ $mS/m$.}
%\label{Hr^(1)}\center
%{\includegraphics[width=0.8\textwidth]{Hr^(1).png}} \quad {%
%\includegraphics[width=0.8\textwidth]{Hr_no_oscillazioni.png}}
%\caption{The imaginary part of the function $g_1(\protect\lambda)$, in the
%case of a $3$-layered underground model with $r=2$ $m$, $f=10$ $kHz$, $m=1$ $A/m^2$ $\protect\sigma_1=33$ $mS/m$, $\protect%
%\sigma_2=20$ $mS/m$, $\protect\sigma_3=100$ $mS/m$, $h_1=2.5$ $m$, $h_2=0.5$ $m$.}
%\label{Hr_no_oscillazioni}
%\end{figure}

\begin{figure}[t]
\centering
%\subfloat[][\label{a}]
%{\includegraphics[width=0.95\textwidth]{Figure_2.eps}} \quad
%\subfloat[][]
%{\includegraphics[width=0.95\textwidth]{Hr_no_oscillazioni.eps}}
\includegraphics[scale=0.5]{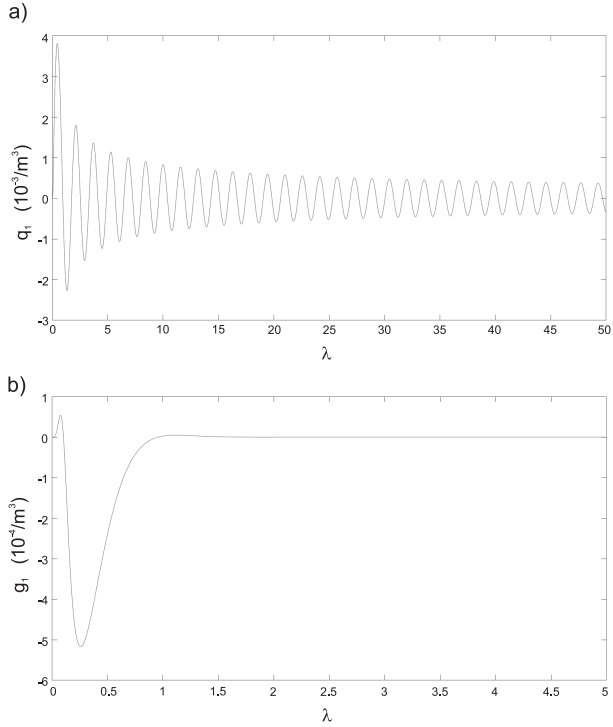}
\caption{(a) The imaginary part of the function $q_1(\lambda)$ and (b) the imaginary part of the function $g_1(\protect\lambda)$, in the case of a $3$-layered underground model with $r=2$ $m$, $f=10$ $kHz$, $m=1$ $A/m^2$ $\protect\sigma_1=33$ $mS/m$, $\protect%
\sigma_2=20$ $mS/m$, $\protect\sigma_3=100$ $mS/m$, $h_1=2.5$ $m$, $h_2=0.5$ $m$.}
\label{Hr_no_oscillazioni}
\end{figure}

Just to provide an example, in Figure \ref{Hr_no_oscillazioni} we plot the imaginary part of $g_1(\lambda)$ and $q_1(\lambda)$ for a given model. It is clear that the oscillations are only retained by the term $q_1(\lambda)$. 

This situation holds true in general and therefore, for suitable positive scalars $s_{l}$, $l=0,1$, we can consider the
following approximations 
\begin{align}
H_{z}^{(N)} & \approx \frac{m}{4\pi }\int_{0}^{s_{0}}g_{0}\left( \lambda \right)
d\lambda +H_{z}^{(1)}, \label{Hz_bis} \\
H_{\rho}^{(N)} & \approx -\frac{m}{4\pi }\int_{0}^{s_{1}}g_{1}\left( \lambda \right)
d\lambda +H_{\rho}^{(1)}, \label{Hr_bis}
\end{align}%
in which we neglect the tail of the integrals.

Theoretically, the truncation error can be bounded as follows.

\begin{proposition} \label{proposizione}
For $l=0,1$ there exists a constant $c$ such that, for $s_l$ large enough,
\begin{equation}
\Big | \int_{s_l}^{\infty} g_l(\lambda) d \lambda \Big | \leq c \sqrt{\frac{1}{8 \pi r}}  \sum_{k=1}^{N-1} | k_{k+1}^2 - k_{k}^2 | \tilde{\gamma}_{k}^{-1} e^{-\tilde{\gamma}_{k} s_l} s_l^{-\frac{1}{2}}.
\end{equation}
\end{proposition}
In order to prove Proposition \ref{proposizione} we need the following lemma.

\begin{lemma} \label{lemma2}
For $u \rightarrow + \infty$,
\begin{equation*}
\int_u^{\infty} e^{-x}x^{-\nu} dx = u^{-\nu}e^{-u} \left( 1+ \mathcal{O} \left( \frac{1}{u} \right) \right).
\end{equation*}
\end{lemma}
\begin{proof}
By using \cite[pag.318]{GR} and \cite[pag.505-504]{Abramowitz} respectively, we have that
\begin{equation*}
\begin{split}
\int_u^{\infty} \frac{e^{-x}}{x^{\nu}} dx &= u^{-\frac{\nu}{2}}e^{-\frac{u}{2}} \mathcal{W}_{-\frac{\nu}{2},\frac{1-\nu}{2}} (u) \\
&=u^{1-\nu}e^{-u} \mathcal{U}(1,2-\nu,u) \\
&= u^{-\nu} e^{-u} \left( 1+ \mathcal{O} \left( \frac{1}{u} \right) \right),
\end{split}
\end{equation*}
where $\mathcal{W}$ is the Whittaker function and $\mathcal{U}$ is the Kummer's confluent hypergeometric function.
\end{proof}

\noindent Now we can prove Proposition \ref{proposizione}.

\begin{proof}
By (\ref{pallino}) for $l=0,1$, we have
\begin{equation*}
\begin{split}
\Big |  \int_{s_l}^{\infty} g_l(\lambda) d \lambda  \Big | \lesssim \sum_{k=1}^{N-1} \frac{| k_{k+1}^2-k_{k}^2 |}{4} \int_{s_l}^{\infty} e^{-\tilde{\gamma}_{k}\lambda}  |J_l(\lambda r)| d \lambda
\end{split}
\end{equation*}
Using the relation (see \cite[pag.364]{Abramowitz})
\begin{equation*}
J_l(t) = \sqrt{\frac{2}{\pi t}} \left[cos \left(t-\frac{1}{2} l \pi-\frac{1}{4} \pi \right) + \mathcal{O} \left(\frac{1}{t} \right) \right], \quad {\rm for} \; t \rightarrow + \infty,
\end{equation*}
and Lemma \ref{lemma2}, we obtain
\begin{equation*}
\begin{split}
\int_{s_l}^{\infty} e^{-\tilde{\gamma}_{k}\lambda}  |J_l(\lambda r)| d \lambda &\lesssim \sqrt{\frac{2}{\pi r}} \int_{s_l}^{\infty} e^{-\tilde{\gamma}_{k}\lambda} \lambda^{-\frac{1}{2}} d \lambda \\
&= \sqrt{\frac{2}{ \pi r}} \int_{\tilde{\gamma}_{k}s_l}^{\infty} e^{-t} t^{-\frac{1}{2}}  dt \\
&\approx \sqrt{\frac{2}{\pi r}} \tilde{\gamma}_{k}^{-1} e^{-\tilde{\gamma}_{k} s_l} s_l^{-\frac{1}{2}},
\end{split}
\end{equation*}
%\begin{equation*}
%\begin{split}
%\int_{s_l}^{\infty} |R_0-\Psi_1| \lambda^2 |J_l( \lambda r)| d \lambda &\approx \int_{s_l}^{\infty} |R_1| \lambda^2 |J_l ( \lambda r)| d \lambda \\ 
%&\approx \sum_{k=1}^{N-1} \frac{|k_{k+1}^2-k_{k}^2|}{4} \int_{s_l}^{\infty} e^{-\tilde{\gamma}_{k}\lambda}  |J_l (\lambda r)| d \lambda \\ 
%&= \sum_{k=1}^{N-1} \frac{|k_{k+1}^2-k_{k}^2|}{4} \int_{s_l}^{\infty} e^{-\tilde{\gamma}_{k}\lambda}  |J_l(\lambda r)| d \lambda \\
%&\lesssim \sqrt{\frac{2}{\pi r}} \sum_{k=1}^{N-1} \frac{|k_{k+1}^2-k_{k}^2|}{4} \int_{s_l}^{\infty} e^{-\tilde{\gamma}_{k}\lambda} \lambda^{-\frac{1}{2}} d \lambda \\
%&= \sqrt{\frac{1}{8 \pi r}} \sum_{k=1}^{N-1} |k_{k+1}^2-k_{k}^2| \tilde{\gamma}_{k}^{-\frac{1}{2}} \int_{\tilde{\gamma}_{k}s_l}^{\infty} e^{-t} t^{-\frac{1}{2}}  dt \\
%&\approx \sqrt{\frac{1}{8 \pi r}} \sum_{k=1}^{N-1} |k_{k+1}^2-k_{k}^2| \tilde{\gamma}_{k}^{-1} e^{-\tilde{\gamma}_{k} s_l} s_l^{-\frac{1}{2}}, 
%\end{split}
%\end{equation*} 
where we have used the symbols $\approx$ and $\lesssim$ to neglect the factor $1+\mathcal{O} \left( \frac{1}{\lambda} \right)$.
\end{proof}

In practice, $s_{l}$ can be taken relatively small to obtain reliable results by a traditional quadrature formula on finite intervals, e.g., $s_{l} = 2\div 3$, for $l=0,1$.  Table \ref{tabella_integrali} shows, for example, the results obtained with $r=2$ $m$ and $f=10$ $kHz$, in the case of a subsurface model composed by 3 layers: $\protect\sigma_1=333$ $mS/m$, $\protect\sigma_2=20$ $mS/m$, $\protect\sigma_3=100$ $mS/m$, $h_1=2.5$ $m$, $h_2=0.5$ $m$.
In Table \ref{tabella_integrali}, as well as in the sequel of this work, for the integral evaluation we have used the Gauss-Kronrod quadrature technique \cite{quadgk}, with relative error tolerance equal to $10^{-8}$.
By adopting the above model, we also compared the fields $H_{z}^{(3)}$ and $H_{\rho}^{(3)}$, computed from the integral formulations (\ref{Hz_bis}) and (\ref{Hr_bis}), with the results of the digital filtering algorithm provided by Singh and Mogi \cite{Emdpler}. Using 300 points, $s_{0} = s_{1} =$ 3, and a tolerance equal to $10^{-9}$ for the integral evaluations, the two curves differ by less than $10^{-8}$. Using more points the error further decreases, but with higher computational costs.

\begin{table}[tbp]
\begin{equation*}
\begin{array}{ccc}
\toprule
s_l & \int_0^{s_0} g_0 (\lambda) d \lambda & \int_0^{s_1} g_1 (\lambda) d \lambda \\
  (m^{-1})      & (A/m) & (A/m) \\
\midrule
0.5 & -0.0871 \cdot 10^{-3} & 0.2797 \cdot 10^{-3}  \\
1.0 & -0.1279 \cdot 10^{-3} & 0.3279 \cdot 10^{-3} \\
1.5 & -0.1316 \cdot 10^{-3} & 0.3283 \cdot 10^{-3} \\ 
2.0 & -0.1317\cdot 10^{-3} & 0.3281 \cdot 10^{-3} \\
2.5 & -0.1317 \cdot 10^{-3} & 0.3280 \cdot 10^{-3} \\
3.0 & -0.1317 \cdot 10^{-3} & 0.3280 \cdot 10^{-3} \\
3.5 & -0.1317 \cdot 10^{-3} & 0.3280 \cdot 10^{-3} \\
4.0 & -0.1317 \cdot 10^{-3} & 0.3280 \cdot 10^{-3} \\
\bottomrule
\end{array}
\end{equation*}%
\caption{Results of the Gauss-Kronrod quadrature for different values of $s_l
$, with $r=2$ $m$ and $f=10$ $kHz$, in the case of a subsurface model composed by 3 layers: $\protect\sigma_1=333$ $mS/m$, $\protect\sigma_2=20$ $mS/m$, $\protect\sigma_3=100$ $mS/m$, $h_1=2.5$ $m$, $h_2=0.5$ $m$.}
\label{tabella_integrali}
\end{table}

\section{Approximation of the fields}

In this section, working with $N=2$ and $N=3$, we derive useful analytical
approximations for $H_{z}^{(N)}$ and $H_{\rho }^{(N)}$, expressed by equations (\ref{eq9}) and (\ref{eq10}). While these approximations can be used in general to have an idea of the main features of the fields, our basic aim is to use them for the inversion problem. First, from (\ref{diff}) we have that
\begin{equation*}
g_{l}\left( \lambda \right) =\frac{4R_{1}u_{1}}{R_{1}k_{1}^{2}+(\lambda
+u_{1})^{2}}\lambda ^{3}J_{l}(\lambda r),\quad l=0,1. 
\end{equation*}%
%Therefore, from the previous splitting
%\begin{eqnarray}
%H_{z}^{(N)} &=&\frac{m}{4\pi }\int_{0}^{\infty }\frac{4R_{1}u_{1}}{%
%R_{1}k_{1}^{2}+(\lambda +u_{1})^{2}}\lambda ^{3}J_{0}(\lambda r)d\lambda
%+H_{z}^{(1)}  \label{Hz_bis} \\
%H_{\rho }^{(N)} &=&-\frac{m}{4\pi }\int_{0}^{\infty }\frac{4R_{1}u_{1}}{%
%R_{1}k_{1}^{2}+(\lambda +u_{1})^{2}}\lambda ^{3}J_{1}(\lambda r)d\lambda
%+H_{\rho }^{(1)}.  \label{Hr_bis}
%\end{eqnarray}%
Before going on, we observe that for $0\leq i,j,\ell \leq N$ 
\begin{eqnarray}
(u_{i}-u_{j})e^{-2u_{\ell}h_{\ell}} &=&\frac{k_{j}^{2}-k_{i}^{2}}{\sqrt{\lambda
^{2}-k_{i}^{2}}+\sqrt{\lambda ^{2}-k_{j}^{2}}}e^{-2\sqrt{\lambda
^{2}-k_{\ell}^{2}}h_{\ell}}  \notag \\
&=&\mathcal{O}\left( \frac{e^{-2\lambda h_{\ell}}}{\lambda }\right) \quad \text{%
for\quad }\lambda \rightarrow \infty   \label{ga}
\end{eqnarray}%
and the same holds for $(\lambda -u_{j})e^{-2u_{\ell}h_{\ell}}$. Actually,
this asymptotic behavior is observed already for $\lambda $ relatively
small, because also the quantities $k_{i}$ are very small. In order to derive approximations of $H_z^{(N)}$ and $H_{\rho}^{(N)}$, in what follows,
whenever possible we neglect the terms involving these factors.

\subsection{The case $N=2$}

For $N=2$ we have%
\begin{equation*}
R_{2}=0,\quad R_{1}=\Psi _{2}e^{-2u_{1}h_{1}}=\frac{u_{1}-u_{2}}{u_{1}+u_{2}}%
e^{-2u_{1}h_{1}},
\end{equation*}%
and therefore%
\begin{equation*}
\frac{4R_{1}u_{1}}{R_{1}k_{1}^{2}+(\lambda +u_{1})^{2}}=\frac{%
4u_{1}(u_{1}-u_{2})}{(\lambda
^{2}-u_{1}^{2})(u_{1}-u_{2})e^{-2u_{1}h_{1}}+(\lambda
+u_{1})^{2}(u_{1}+u_{2})}e^{-2u_{1}h_{1}}.
\end{equation*}%
By (\ref{ga}) we can cosider the approximation%
\begin{equation*}
(\lambda ^{2}-u_{1}^{2})(u_{1}-u_{2})e^{-2u_{1}h_{1}}\approx 0.
\end{equation*}%
Moreover, using the first order approximation%
\begin{equation}
u_{1}=\sqrt{\lambda ^{2}+i\omega \mu \sigma _{1}}=\lambda \sqrt{1+\frac{%
i\omega \mu \sigma _{1}}{\lambda ^{2}}}\approx \lambda \left( 1+\frac{%
i\omega \mu \sigma _{1}}{2\lambda ^{2}}\right) , \label{approssimazione}
\end{equation}%
we obtain 
\begin{align}
\frac{4R_{1}u_{1}}{R_{1}k_{1}^{2}+(\lambda +u_{1})^{2}}& \approx \frac{%
4u_{1}(u_{1}-u_{2})}{(\lambda +u_{1})^{2}(u_{1}+u_{2})}e^{-2u_{1}h_{1}}
\label{9} \\
& =\frac{4u_{1}(u_{1}-u_{2})^{2}(\lambda -u_{1})^{2}}{(\lambda
^{2}-u_{1}^{2})^{2}(u_{1}^{2}-u_{2}^{2})}e^{-2u_{1}h_{1}} \\
& \approx \frac{i\omega \mu (\sigma _{1}-\sigma _{2})}{4\lambda ^{3}}\left(
1+\frac{i\omega \mu \sigma _{1}}{2\lambda ^{2}}\right) e^{-2u_{1}h_{1}}.
\label{asterisco}
\end{align}

At this point, we use the non standard approximation 
\begin{equation}
u_{1}=\sqrt{\lambda ^{2}+i\omega \mu \sigma _{1}}\approx \lambda +\sqrt{%
i\omega \mu \sigma _{1}}.  \label{approx}
\end{equation}%
The main reason for this choice is that the standard (\ref{approssimazione}) leads to final approximation that does not allow to simplify the integrals with the existing formulas. Anyway the approximation (\ref{approx}) is partially justified by observing that 
\begin{equation*}
\frac{u_{1}}{\lambda +\sqrt{i\omega \mu \sigma _{1}}}\rightarrow 1,\quad 
\text{for\quad }\lambda \rightarrow 0\text{ and }\lambda \rightarrow +\infty.
\end{equation*}%
Now, since%
\begin{equation*}
\lambda +\sqrt{i\omega \mu \sigma _{1}}=\lambda +\frac{\sqrt{2}}{2}\sqrt{%
\omega \mu \sigma _{1}}+i\frac{\sqrt{2}}{2}\sqrt{\omega \mu \sigma _{1}},
\end{equation*}%
we obtain 
\begin{eqnarray*}
e^{-2u_{1}h_{1}} &\approx &e^{-2\lambda h_{1}-h_{1}\sqrt{2\omega \mu \sigma
_{1}}}[\cos (h_{1}\sqrt{2\omega \mu \sigma _{1}})-i\sin (h_{1}\sqrt{2\omega
\mu \sigma _{1}})] \\
&\approx &e^{-2\lambda h_{1}-h_{1}\sqrt{2\omega \mu \sigma _{1}}}.
\end{eqnarray*}

Using the above approximation in (\ref{asterisco}) we have 
\begin{equation*}
\Im \left( \frac{4R_{1}u_{1}}{R_{1}k_{1}^{2}+(\lambda +u_{1})^{2}}\right)
\approx \frac{\omega \mu (\sigma _{1}-\sigma _{2})e^{-2\lambda h_{1}-h_{1}%
\sqrt{2\omega \mu \sigma _{1}}}}{4\lambda ^{3}}.
\end{equation*}%
where the symbol $\Im$ indicates the imaginary part. Moreover, by the general formula (see \cite{GR}, p. $707$)%
\begin{eqnarray*}
\int_{0}^{\infty }e^{-\alpha x}J_{\nu }(\beta x)dx &=&\frac{\beta ^{-\nu }[%
\sqrt{\alpha ^{2}+\beta ^{2}}-\alpha ]^{\nu }}{\sqrt{\alpha ^{2}+\beta ^{2}}}%
,\quad \Re (\nu )>-1,\Re (\alpha \pm i\beta )>0 ,
\end{eqnarray*}%
where the symbol $\Re$ indicates the real part. We finally obtain, for the integral in equation (\ref{eq9}), the following approximation:
\begin{align*}
\Im \left( \int_{0}^{\infty }\frac{4R_{1}u_{1}}{R_{1}k_{1}^{2}+(\lambda
+u_{1})^{2}}\lambda ^{3}J_{0}(\lambda r)d\lambda \right) & \approx \frac{%
\omega \mu (\sigma _{1}-\sigma _{2})e^{-h_{1}\sqrt{2\omega \mu \sigma _{1}}}%
}{4}\int_{0}^{\infty }e^{-2\lambda h_{1}}J_{0}(\lambda r)d\lambda \\
& =\frac{\omega \mu (\sigma _{1}-\sigma _{2})e^{-h_{1}\sqrt{2\omega \mu
\sigma _{1}}}}{4\sqrt{4h_{1}^{2}+r^{2}}}.
\end{align*}

Using the same arguments for $H_{\rho}^{(2)}$, the integral in equation (\ref{eq10}) can be approximated as 
\begin{align*}
&\Im \left( \int_{0}^{\infty }\frac{-4R_{1}u_{1}}{R_{1}k_{1}^{2}+(\lambda
+u_{1})^{2}}\lambda ^{3}J_{1}(\lambda r)d\lambda \right)  \approx \\
&\approx -\frac{%
\omega \mu (\sigma _{1}-\sigma _{2})e^{-h_{1}\sqrt{2\omega \mu \sigma _{1}}}%
}{4}\int_{0}^{\infty }e^{-2\lambda h_{1}}J_{1}(\lambda r)d\lambda \\
& =-\frac{\omega \mu (\sigma _{1}-\sigma _{2})e^{-h_{1}\sqrt{2\omega \mu
\sigma _{1}}}(\sqrt{4h_{1}^{2}+r^{2}}-2h_{1})}{4r\sqrt{4h_{1}^{2}+r^{2}}}.
\end{align*}

Finally, 
\begin{equation*}
\Im (H_{z}^{(2)}) \approx \frac{m}{4\pi }\left( \frac{\omega \mu (\sigma
_{1}-\sigma _{2})e^{-h_{1}\sqrt{2\omega \mu \sigma _{1}}}}{4\sqrt{%
4h_{1}^{2}+r^{2}}}\right) +\Im \left( H_{z}^{(1)}\right)
\end{equation*}
and
\begin{equation*}
\Im (H_{\rho}^{(2)}) \approx -\frac{m}{4\pi }\left( \frac{\omega \mu
(\sigma _{1}-\sigma _{2})e^{-h_{1}\sqrt{2\omega \mu \sigma _{1}}}(\sqrt{%
4h_{1}^{2}+r^{2}}-2h_{1})}{4r\sqrt{4h_{1}^{2}+r^{2}}}\right) +\Im \left(
H_{\rho}^{(1)}\right).
\end{equation*}

\subsection{The case $N=3$}

For $N=3$ we have 
\begin{align*}
R_{3}& =0, \\
R_{2}& =\Psi _{3}e^{-2u_{2}h_{2}}=\frac{u_{2}-u_{3}}{u_{2}+u_{3}}%
e^{-2u_{2}h_{2}}, \\
R_{1}& =\frac{\frac{u_{2}-u_{3}}{u_{2}+u_{3}}e^{-2u_{2}h_{2}}+\frac{%
u_{1}-u_{2}}{u_{1}+u_{2}}}{\frac{u_{2}-u_{3}}{u_{2}+u_{3}} \frac{%
u_{1}-u_{2}}{u_{1}+u_{2}}e^{-2u_{2}h_{2}}+1} e^{-2u_{1}h_{1}} \\
& =\frac{%
(u_{1}+u_{2})(u_{2}-u_{3})e^{-2u_{2}h_{2}}+(u_{1}-u_{2})(u_{2}+u_{3})}{%
(u_{2}-u_{3})(u_{1}-u_{2})e^{-2u_{2}h_{2}}+(u_{1}+u_{2})(u_{2}+u_{3})}
e^{-2u_{1}h_{1}} \\
& \approx \frac{%
(u_{1}+u_{2})(u_{2}-u_{3})e^{-2u_{2}h_{2}}+(u_{1}-u_{2})(u_{2}+u_{3})}{%
(u_{1}+u_{2})(u_{2}+u_{3})} e^{-2u_{1}h_{1}},
\end{align*}%
where, as before, we have used (\ref{ga}). We obtain%
\begin{eqnarray}
\frac{4R_{1}u_{1}}{R_{1}k_{1}^{2}+(\lambda +u_{1})^{2}}
&=&4u_{1}[(u_{1}+u_{2})(u_{2}-u_{3})e^{-2(u_{1}h_{1}+u_{2}h_{2})}+(u_{1}-u_{2})(u_{2}+u_{3})e^{-2u_{1}h_{1}}]\times
\notag \\
&&[(u_{1}+u_{2})(u_{2}-u_{3})(\lambda
^{2}-u_{1}^{2})e^{-2(u_{1}h_{1}+u_{2}h_{2})}+(u_{1}-u_{2})(u_{2}+u_{3})(%
\lambda ^{2}-u_{1}^{2})e^{-2u_{1}h_{1}}  \notag \\
&&+(\lambda+u_1)^2(u_1-u_2)(u_2-u_3)e^{-2u_2h_2}+(\lambda +u_{1})^{2}(u_{1}+u_{2})(u_{2}+u_{3})]^{-1}  \notag \\
&\approx &4u_{1}\frac{%
(u_{1}+u_{2})(u_{2}-u_{3})e^{-2(u_{1}h_{1}+u_{2}h_{2})}+(u_{1}-u_{2})(u_{2}+u_{3})e^{-2u_{1}h_{1}}%
}{(\lambda +u_{1})^{2}(u_{1}+u_{2})(u_{2}+u_{3})}  \label{4} \\
&=&\frac{4u_{1}(u_{2}-u_{3})e^{-2(u_{1}h_{1}+u_{2}h_{2})}}{(\lambda
+u_{1})^{2}(u_{2}+u_{3})}+\frac{4u_{1}(u_{1}-u_{2})e^{-2u_{1}h_{1}}}{%
(\lambda +u_{1})^{2}(u_{1}+u_{2})},  \label{17}
\end{eqnarray}%
where the approximation (\ref{4}) arises again from (\ref{ga}). We observe
that the two terms of (\ref{17}) are very similar and, in particular, the
second one corresponds to (\ref{9}) for $N=2$. Therefore, using the
approximations (\ref{approssimazione}) and (\ref{approx}) for $u_{1}$, $%
u_{2} $ and $u_{3}$, we obtain  
\begin{align*}
&\Im \left( \int_{0}^{\infty} \frac{4R_{1}u_{1}}{R_{1}k_{1}^{2}+(\lambda+u_{1})^{2}} \lambda ^{3}J_{0}(\lambda r)d\lambda \right) \approx \\
& \approx \frac{\omega \mu e^{-h_1 \sqrt{2 \omega \mu \sigma_1}}}{4} \Big ( \frac{(\sigma_{2}-\sigma _{3})e^{-h_{2}\sqrt{2\omega \mu \sigma _{2}}}}{\sqrt{4(h_{1}+h_{2})^{2}+r^{2}}} +\frac{\sigma _{1}-\sigma _{2}}{\sqrt{4h_{1}^{2}+r^{2}}} \Big ).
\end{align*}

Analougsly,

\begin{eqnarray*}  
&&\Im \left( \int_{0}^{\infty }\frac{-4R_{1}u_{1}}{R_{1}k_{1}^{2}+(\lambda
+u_{1})^{2}}\lambda ^{3}J_{1}(\lambda r)d\lambda \right) \approx \\ &\approx& \frac{%
\omega \mu e^{-h_{1}\sqrt{2\omega \mu \sigma _{1}}}}{-4r}\left( \frac{%
(\sigma _{2}-\sigma _{3})e^{-h_{2}\sqrt{2\omega \mu \sigma _{2}}}[\sqrt{%
4(h_{1}+h_{2})^{2}+r^{2}}-2(h_{1}+h_{2})]}{\sqrt{4(h_{1}+h_{2})^{2}+r^{2}}}%
\right. \\
&&\left. +\frac{(\sqrt{4h_{1}^{2}+r^{2}}-2h_{1})(\sigma _{1}-\sigma _{2})}{%
\sqrt{4h_{1}^{2}+r^{2}}}\right).
\end{eqnarray*}

%\begin{table}[tbp]
%\begin{equation*}
%\begin{array}{ccc}
%\toprule
%r & E_z & E_{{\rho}} \\
%(m) & (\%) & (\%) \\
%\midrule
%1 & 0.0 & 0.2 \\
%2 & 0.0 & 0.8 \\
%3 & 0.0 & 1.6 \\ 
%4 & 0.0 & 2.4 \\
%5 & 0.0 & 3.1 \\
%6 & 0.0 & 3.5 \\
%7 & 0.0 & 3.7 \\
%8 & 0.0 & 3.6 \\
%9 & 0.1 & 3.1 \\
%10 & 0.1 & 2.3 \\
%\bottomrule
%\end{array}
%\end{equation*}%
%\caption{Relative percentage errors between $\Im(H_{\rho}^{(3)})$ and $\Im(H_z^{(3)})$ and their approximations in (\ref{puntino}) and (\ref{ImHr}), in the case of a $3$-layered underground model with $f=10$ $kHz$, $\sigma_1=50$ $mS/m$, $\sigma_2=4.9$ $mS/m$, $\sigma_3=18.2$ $mS/m$, $h_1=2.5$ $m$, $h_2=0.5$ $m$.}
%\label{tabella_errori}
%\end{table}

%\begin{table}[tbp]
%\begin{equation*}
%\begin{array}{ccccccccccc}
%\toprule
%r \; (m) & 1 & 2 & 3 & 4 & 5 & 6 & 7 & 8 & 9 & 10 \\
%\midrule
%E_{z} \; (\%) & 0.0 & 0.0 & 0.0 & 0.0 & 0.0 & 0.0 & 0.0 & 0.0 & 0.1 & 0.1 \\
%E_{\rho } \; (\%) & 0.2 & 0.8 & 1.6 & 2.4 & 3.1 & 3.5 & 3.7 &
%3.6 & 3.1 & 2.3 \\
%\bottomrule
%\end{array}
%\end{equation*}%
%\caption{Relative percentage errors between $\Im(H_{\rho}^{(3)})$ and $\Im(H_z^{(3)})$ and their approximations in (\ref{puntino}) and (\ref{ImHr}), in the case of a $3$-layered underground model with $f=10$ $kHz$, $\sigma_1=50$ $mS/m$, $\sigma_2=4.9$ $mS/m$, $\sigma_3=18.2$ $mS/m$, $h_1=2.5$ $m$, $h_2=0.5$ $m$.}
%\label{tabella_errori}
%\end{table}

Finally,
\begin{equation} \label{puntino}
\Im (H_{z}^{(3)})\approx \frac{m}{4\pi } \frac{\omega \mu
e^{-h_{1}\sqrt{2\omega \mu \sigma _{1}}}}{4} \left( \frac{(\sigma
_{2}-\sigma _{3})e^{-h_{2}\sqrt{2\omega \mu \sigma _{2}}}}{\sqrt{%
4(h_{1}+h_{2})^{2}+r^{2}}}+\frac{\sigma _{1}-\sigma _{2}}{\sqrt{%
4h_{1}^{2}+r^{2}}}\right) +\Im \left( H_{z}^{(1)}\right)
\end{equation}%
and 
\begin{equation} \label{ImHr}
\begin{split}
\Im (H_{\rho}^{(3)}) &\approx -\frac{m}{4\pi }\frac{\omega \mu e^{-h_{1} \sqrt{2\omega \mu \sigma _{1}}}}{4r} \Big( \frac{(\sigma _{2}-\sigma
_{3})e^{-h_{2}\sqrt{2\omega \mu \sigma _{2}}}[\sqrt{4(h_{1}+h_{2})^{2}+r^{2}}-2(h_{1}+h_{2})]}{\sqrt{4(h_{1}+h_{2})^{2}+r^{2}}} \\
& +\frac{(\sqrt{4h_{1}^{2}+r^{2}}-2h_{1})(\sigma _{1}-\sigma _{2})}{\sqrt{4h_{1}^{2}+r^{2}}} \Big) +\Im (H_{\rho}^{(1)}).
\end{split} 
\end{equation}

%In Table \ref{tabella_errori} we show the relative error (percentage) between the imaginary part of (\ref{Hz_bis}) and (\ref{Hr_bis}), with $N=3$, and their approximations  (\ref{puntino}) and (\ref{ImHr}), for a given set of parameters and different values of the offset.

\section{The inverse problem}

Let $p=(\sigma _{1},...,\sigma _{N},h_{1},...,h_{N-1})\in \mathbb{R}^{2N-1}$
be the vector of underground parameters. The components of the magnetic field at the surface are functions of these parameters and the offset $r$, thus $%
H_{z}^{(N)}=H_{z}^{(N)}(p,r)$ and $H_{\rho}^{(N)}=H_{\rho}^{(N)}(p,r)$.
Given a vector of observations $d=\left( d_{z,1},...,d_{z,k},d_{\rho,1},...,d_{\rho,k}\right) ^{T}$, corresponding to distances $r_{1},...,r_{k}$ ($2k>2N-1$), the inverse
problem can be formulated as%
\begin{equation}
\min_{p}\sum_{i=1}^{k}\left\{ \left[ \Im (H_{z}^{(N)}(p,r_{i}))-d_{z,i}%
\right] ^{2}+\left[ \Im (H_{\rho}^{(N)}(p,r_{i}))-d_{\rho,i}\right]
^{2}\right\}.  \label{minp}
\end{equation}%
It is known that other kind of minimizations are possible, that is, by using a
different norm, and eventually a regularization term can be included to
reduce the effect of noise on the measuraments. By defining $\mathcal{H}:\mathbb{R}%
^{2N-1}\rightarrow \mathbb{R}^{2k}$ as%
\begin{equation} \label{H_cal}
\mathcal{H}(p)=\left( \Im (H_{z}^{(N)}(p,r_{1})),...,\Im (H_{z}^{(N)}(p,r_{k})),\Im
(H_{\rho}^{(N)}(p,r_{1})),...,\Im (H_{\rho}^{(N)}(p,r_{k}))\right) ^{T},
\end{equation}%
the problem (\ref{minp}) can be rewritten in the compact form%
\begin{equation}
\min_{p}\left\Vert \mathcal{H}(p)-d\right\Vert ^{2}.  \label{minc}
\end{equation}%
Due to the complexity of the function $\mathcal{H}(p)$, whose computation requires the evaluation of integrals, it is necessary to employ a derivative free minimization algorithms based, for example, on a BFGS line-search method (\cite{Broyden,Fletcher,Goldfarb,Shanno}) and on the SA global-search technique (\cite{Kirkpatrick,Goffe}). However, these methods are very demanding in terms of computational resources and elaboration time. By using the approximations given in the previous section, below we present a more efficient algorithm to solve (\ref%
{minc}) in the nontrivial case of $N=3$.

%\begin{align*}
%\Im (H_{z}) &= \frac{m}{4\pi }F_{z}(p)+ \Im (H_{z}^{(1)}) \approx \frac{m}{4\pi }%
%\bar{F}_{z}(p)+ \Im (H_{z}^{(1)}), \\
%\Im (H_{\rho}) &= \frac{m}{4\pi }F_{\rho}(p)+ \Im (H_{\rho}^{(1)}) \approx \frac{m}{%
%4\pi }\bar{F}_{\rho}(p)+ \Im (H_{\rho}^{(1)}),
%\end{align*}
%where 
%\begin{align*}
%F_z(p)&=\Im \left( \int_{0}^{\infty } \frac{4R_{1}u_{1}}{R_{1}k_{1}^{2}+(%
%\lambda +u_{1})^{2}}\lambda ^{3}J_{0}(\lambda r)d\lambda \right) , \\
%F_{\rho}(p)&=\Im \left( \int_{0}^{\infty }-\frac{4R_{1}u_{1}}{%
%R_{1}k_{1}^{2}+(\lambda +u_{1})^{2}}\lambda ^{3}J_{1}(\lambda r)d\lambda
%\right) ,
%\end{align*}
%and 
%\begin{align*}
%\bar{F}_z(p)&= \left( \frac{\omega \mu e^{-h_{1}\sqrt{2\omega \mu \sigma _{1}%
%}}}{4}\right) \frac{(\sigma _{2}-\sigma _{3})e^{-h_{2}\sqrt{2\omega \mu
%\sigma _{2}}}}{\sqrt{4(h_{1}+h_{2})^{2}+r^{2}}}+\frac{\sigma _{1}-\sigma _{2}%
%}{\sqrt{4h_{1}^{2}+r^{2}}}, \\
%\bar{F}_{\rho}(p)&=-\frac{\omega \mu e^{-h_{1}\sqrt{2\omega \mu \sigma _{1}}}%
%}{4r} \Big[ \frac{(\sigma _{2}-\sigma _{3})e^{-h_{2}\sqrt{2\omega \mu \sigma
%_{2}}} [\sqrt{4(h_{1}+h_{2})^{2}+r^{2}}-2(h_{1}+h_{2})]}{\sqrt{%
%4(h_{1}+h_{2})^{2}+r^{2}}} \\
%&+ \frac{(\sqrt{4h_{1}^{2}+r^{2}}-2h_{1})(\sigma _{1}-\sigma _{2})}{ \sqrt{%
%4h_{1}^{2}+r^{2}}}\Big].
%\end{align*}

Let $p=(\sigma _{1},\sigma _{2},\sigma _{3},h_{1},h_{2}) \in \mathbb{R}^{5}$ be the vector of
the underground parameters and consider the fields $H_z^{(3)}=H_z$ and $%
H_{\rho}^{(3)}=H_{\rho}$, defined by equations (\ref{eq9}) and (\ref%
{eq10}). From the previous section, $\Im(H_z)$ and $\Im (H_{\rho})$ can be approximated using equations (\ref{puntino}) and (\ref{ImHr}) at different offsets $r_1, \ldots, r_k$. We indicate these approximations as $L_z$ and $L_{\rho}$.

Suppose a real underground model be characterized by the parameters vector $p^{\star}= (\sigma^{\star}_1, \sigma^{\star}_2, \sigma^{\star}_3, h^{\star}_1, h^{\star}_2) \in \mathbb{R}^{5}$. 
%The observations vector $d$ can be computed either using the integral formulations \ref{Hz_bis} and \ref{Hr_bis}, or the digital filtering algorithm (e.g., \cite{Emdpler}). These forward numerical simulations provide the best representation of realistic EM surveys.
The observations vector $d$ is computed using the integral formulations (\ref{Hz_bis}) and (\ref{Hr_bis}), being numerical simulations the best representation of realistic EM surveys.
Given an initial guess $p_0$ of $p^{\star}$, we search for a certain $\bar{p} \in \mathbb{R}^{5}$ such that 
\begin{equation}  \label{minimo}
\bar{p} = (\bar{\sigma}_1, \bar{\sigma}_2, \bar{\sigma}_3,
\bar{h}_1, \bar{h}_2) = \arg \min_{p}\left\Vert \bar{\mathcal{H}}(p)-d\right\Vert ^{2},
\end{equation}
where $\bar{\mathcal{H}} :\mathbb{R}^{5}\rightarrow \mathbb{R}^{2k}$ is defined as 
\begin{equation*}
\bar{\mathcal{H}}(p)=\left( L_z(p,r_{1})),...,L_z(p,r_{k})),L%
_{\rho}(p,r_{1})),...,L_{\rho}(p,r_{k}))\right) ^{T}.
\end{equation*}
At this point we use $\bar{p}$ as initial guess for a second minimization procedure with the integral formulations (\ref{Hz_bis}) and (\ref{Hr_bis}) instead of the analytical approximations (\ref{puntino}) and (\ref{ImHr}).
Therefore, we look for $\hat{p} \in \mathbb{R}^{5}$ such that 
\begin{equation}
\hat{p}=(\hat{\sigma}_{1},\hat{\sigma}_{2},%
\hat{\sigma}_{3},\hat{h}_{1},\hat{h}_{2})=\arg \min_{p}\left\Vert \mathcal{H}(p)-d\right\Vert ^{2},  \label{secondo_minimo}
\end{equation}
where $\mathcal{H}:\mathbb{R}^{5}\rightarrow \mathbb{R}^{2k}$ is defined by (\ref{H_cal}).
Finally, the approximate solution of the inverse problem is given by the vector $\hat{p} \in \mathbb{R}^{5}$.

The above procedure can be summarized by the following algorithm. 

\begin{algorithm} \label{algoritmo}
Given $p_0 \in \mathbb{R}^5$ and $d \in \mathbb{R}^{2k}$,
\begin{enumerate}
\item solve the problem (\ref{minimo}) to find $\bar{p}$,
\item solve the problem (\ref{secondo_minimo}), with $\bar{p}$ as initial guess, to find $\hat{p}$.
\end{enumerate}
\end{algorithm}

Since $\bar{\mathcal{H}}(p) \approx \mathcal{H}(p)$ is a good approximation (see the examples in the following section), the double step minimization allows to considerably reduce the computational cost because $\bar{\mathcal{H}}$ does not require the evaluation of integrals.

\section{Examples}

In this section we present the results obtained by implementing the Algorithm \ref{algoritmo}, using the BFGS optimization solver. We compare these results with the solutions provided by the SA optimization solver. The examples aim to simulate an EM acquisition on river levees using the DUALEM system (https://dualem.com), for which the offsets are $r_1 = 2$ $m$, $r_2 = 4$ $m$, $r_{3} = 6$ $m$ and $r_{4} = 8$ $m$, and the frequency is about $f = 10$ $kHz$.

\subsection {Conductivity models of river levees}

Regarding the choice of the underground models, we considered a specific application related to the internal composition of river levees, which may collapse due to the condition of the soils that form the embankments. Particularly, the presence of gravel lenses inside the embankment body constitutes a critical factor for the structural levee stability. 
Therefore, the subsoil models are composed by three layers, representing a highly porous gravel level embedded in sediments composed by clay and silty sands.
Two (extreme) cases are considered: the first case (Models 1 and 3) represents a dry levee (summer conditions) and the second case (Models 2 and 4) represents a wet levee (winter conditions). Moreover, different layer thicknesses are taken into account.
The layer conductivities are computed with the complex refractive index model (CRIM), described in \ref{AppA}, and the material properties shown in Table \ref{Sigmamat} and Table \ref{tabella_modelli}.
The computed values in Table \ref{tabella_modelli} represent typical conductivities of shallow sediments, frequently measured on river embankments (e.g., \cite{Zohdy, Francese}).   

%\begin{table}[tbp]
%\begin{equation*}
%\begin{array}{lc}
%\toprule
%     &  \sigma (S/m) \\
%\midrule
%Clay & 0.2 \cite{Carcione2012} \\
%Sand/Silt grains & 0.01 \cite{Carcione2012} \\
%Tap water & 0.1 \cite{Palacky} \\
%Air & 0.0001 \cite{Pawar} \\
%\bottomrule
%\end{array}
%\end{equation*}%
%\caption{Material electrical conductivity}
%\label{Sigmamat}
%\end{table}

\begin{table}[tbp]
\begin{center}
\begin{tabular}{lc}
\toprule
     &  $\sigma (S/m)$ \\
\midrule
Clay & $0.2$ \cite{Carcione2012} \\
Sand/Silt grains & $0.01$ \cite{Carcione2012} \\
Tap water & $0.1$ \cite{Palacky} \\
Air & $0.0001$ \cite{Pawar}  \\
\bottomrule
\end{tabular}
\caption{Material electrical conductivity}
\label{Sigmamat}
\end{center}
\end{table}

\begin{table}[h]
\begin{equation*}
\begin{array}{cccccccc}
\toprule
\text{Medium}  & \text{Layer}  &  \text{$C$ ($\%$)} & \text{$\phi$ ($\%$)} & \text{$S_w$ ($\%$)} &  \text{$\sigma$ ($mS/m$)} &  \text{$h$ ($m$)} & \text{Lithology} \\ 
\midrule
\multirow{3}*{Model 1} & 1 & 50 & 20 & 3 & 50.0 & 2.5 & \text{Dry silt and clay} \\
                       & 2 & 1 & 37 & 1 &  4.9 & 0.5 & \text{Dry gravel lens} \\
                       & 3 & 25 & 30 & 2 & 18.2 &     & \text{Dry sand/silt and clay} \\
\midrule
\multirow{3}*{Model 2} & 1 & 50 & 20 & 92 & 76.9 & 2.5 & \text{Wet silt and clay} \\
                       & 2 & 1 & 37 & 98 & 32.3 & 0.5 & \text{Wet gravel lens} \\
                       & 3 & 25 & 30 & 98 & 50.0 &     & \text{Wet sand/silt and clay} \\
\midrule
\multirow{3}*{Model 3} & 1 & 50 & 20 & 3 & 50.0 & 3.0 & \text{Dry silt and clay} \\
                       & 2 & 1 & 37 & 1 & 4.9  & 2.0 & \text{Dry gravel lens} \\
                       & 3 & 25 & 30 & 2 & 18.2 &     & \text{Dry sand/silt and clay} \\
\midrule
\multirow{3}*{Model 4} & 1 & 50 & 20 & 92 & 76.9 & 3.0 & \text{Wet silt and clay} \\
                       & 2 & 1 & 37 & 98 & 32.3 & 2.0 & \text{Wet gravel lens} \\
                       & 3 & 25 & 30 & 98 & 50.0 &     & \text{Wet sand/silt and clay} \\
\bottomrule &  &  &   
\end{array}
\end{equation*}%
\caption{Petrophysical properties of the river levees models used in the numerical forward simulations.}
\label{tabella_modelli}
\end{table}

\subsection {Results of the EM modelling and inversion}

\begin{figure}[t]
\begin{center}
\includegraphics[scale=0.6]{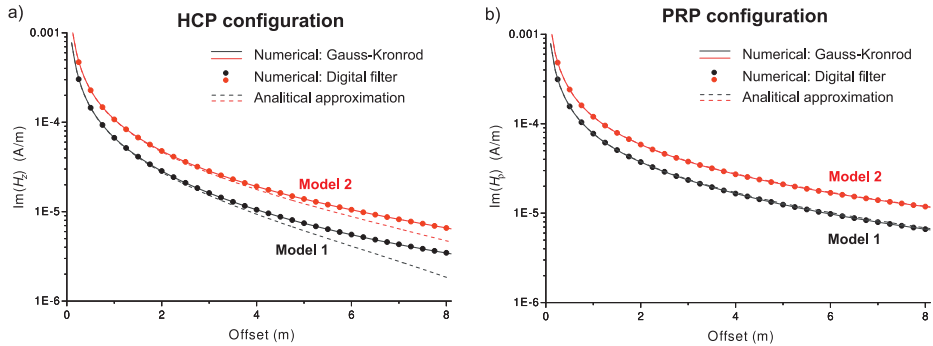}
\end{center}
\caption{Comparison between the fields $H_{z}^{(3)}$ (a) and $H_{\rho}^{(3)}$ (b) computed adopting the numerical Gauss-Kronrod quadrature technique (solid lines), the analytical approximation (dashed lines) and the numerical digital filtering algorithm provided by Singh and Mogi \cite{Emdpler} (symbols). The two components of the the magnetic field correspond to the HCP and PRP coil configurations, respectively. Both Model 1 and Model 2 are considered.}
\label{EM_Curves}
\end{figure}

First of all, we compared the fields $H_{z}^{(N)}$ and $H_{\rho}^{(N)}$ computed adopting the numerical Gauss-Kronrod quadrature technique, the analytical approximation and the numerical digital filtering algorithm provided by Singh and Mogi \cite{Emdpler}. Figure \ref{EM_Curves} shows the imaginary part of (\ref{Hz_bis}) and (\ref{Hr_bis}), with $N=3$, and their analytical approximations  (\ref{puntino}) and (\ref{ImHr}), for both Model 1 and Model 2. 

The two components of the the magnetic field correspond to the HCP and PRP coil configurations, respectively. The values computed using the numerical digital filer approach are shown as well. This plot confirms the very good concordance of the two numerical representations, whose difference was previously estimated by less than $10^{-8}$. Moreover, the curves show that the analytical approximation is better for the PRP configuration.  \\
To solve the minimization problems (\ref{minimo}) and (\ref{secondo_minimo}), we used the quasi-Newton method BFGS, with step tolerance and termination tolerance in the order of the machine precision. We applied the Algorithm \ref{algoritmo} to the observations vector $d$ obtained from the models defined in Table \ref{tabella_modelli} and calculated the relative percentage error between the solution $\hat{p}$ and the real underground model parameters $p^{\star}$. We also applied the SA global-search technique to directly solve the minimization problem (\ref{secondo_minimo}), i.e., without using the analytical approximations (\ref{puntino}) and (\ref{ImHr}). In fact, this kind of technique do not benefit from the two-step procedure described in the previous section and schematized in Algorithm \ref{algoritmo}, because the final result does not depend on the initial model.\\ The observation data vector $d$ is generated by
\begin{equation*}
d=\mathcal{H}(p^{\star}) + \eta,
\end{equation*} 
where $\eta$ is a random white noise vector such that the Noise-to-Signal Ratio ($NSR$) is $\epsilon = \frac{\| \eta \|}{\| d \|}$. We consider two values of percentage $NSR$, $\epsilon = 0.1 \%, 0.5 \%$. These values are compatible with the results of laboratory test, showing that environmental RMS noise levels, in terms of apparent conductivity measurements at 8 m offset, are typically less than 1 mS/m (see https://dualem.com).
The solution errors were averaged over 20 simulations for each model, in order to have an estimate of the parameters expectation values.    

For the SA optimization method, the initial temperature is 10$^{6}$, the temperature reduction factor is 0.1 and the error tolerance for termination is 10$^{-9}$, without noise. By adding the random noise to the data, the error tolerance was increased to 10$^{-6}$. The lower and upper bounds are fixed, on the basis of observations of river levees carried out with other direct and indirect methods, to the following values: 3 $mS/m$ and 1 $S/m$ for the conductivity, and to 0.1 $m$ and 4 $m$ for the thickness. \\ 

The results of the simulations are shown in Tables \ref{tabella_risultati1}-\ref{tabella_risultati2}-\ref{tabella_risultati3}, where we report the mean relative error over $20$ experiments for each model in order to reduce the random dependence on the data. The optimization solvers always converge for all the simulations. The errors are mainly due to the presence of a large number of local minima in the objective function. The simulations without random noise (Table \ref{tabella_risultati1}) show the best results, with the maximum average error around 5 \%. Then, as expected, the errors increase with increasing $NSR$, reaching an average maximum error of $\approx 35 \%$ for $NSR = 0.5$ (Table \ref{tabella_risultati3}). The conductivity and thickness of the intermediate layer are affected by the larger errors. The results of SA show very low errors on the top and bottom layers, while the errors related to the BFGS method are more homogeneously distributed on the three layers. Moreover, thinner layers show larger errors.  \\ 

Overall, the performance of the two methods are very similar in terms of accuracy (see average errors at the bottom of Tables \ref{tabella_risultati1}-\ref{tabella_risultati2}-\ref{tabella_risultati3}). In the absence of noise, the SA technique provides slightly lower errors than BFGS, while in presence of noise we observe the opposite behavior. However, it is very important to note that the Algorithm \ref{algoritmo} with the BFGS optimization solver is more than 10 times faster than SA, requiring less than 1 s for a single inversion on a common laptop.    

%\begin{table}[tbp]
%\begin{equation*}
%\begin{array}{ccccccc}
%\toprule
% & \multicolumn{2}{c}{\mathrm{\textbf{Underground} \; \textbf{model}}} & \multicolumn{2}{c}{\mathrm{\textbf{Mean} \; \textbf{relative} \; \textbf{error}} (6m)} & \multicolumn{2}{c}{\mathrm{\textbf{Mean} \; \textbf{relative} \; \textbf{error}} (8m)}  \\
% & \mathrm{Conductivity} & \mathrm{Thickness} & \mathrm{Conductivity} & \mathrm{Thickness} & \mathrm{Conductivity} & \mathrm{Thickness} \\
% &(mS/m) & (m) & (mS/m) & (m) & (\%) & (\%) \\ 
%\midrule
%\multirow{3}*{Model 1} &  50.0 & 2.5 & 53.2 & 25.2 & 24.8 & 18.9 \\
%&  4.9 & 0.5 & >100 & 59.7 & 46.0 & 7.7 \\
%& 18.2 & & 49.5 & & 18.4 & \\
%\midrule
%\multirow{3}*{Model 2} &  76.9 & 2.5 & 24.2 & 18.7 & 13.9 & 16.1 \\
%& 32.3 & 0.5 & >100 & 29.4 & 33.2 & 44.4 \\
%& 50.0 & & 16.3 & & 8.6 & \\
%\midrule
%\multirow{3}*{Model 3} &  50.0 & 3.0 & 35.0 & 16.3 & 23.4 & 15.4 \\
%& 4.9 & 2.0 & >100 & 20.4 &  >100 & 21.4 \\
%& 18.2 & & 22.6 & & 9.4 & \\
%\midrule
%\multirow{3}*{Model 4} & 76.9 & 3.0 & 30.7 & 16.7 &  20.4 & 13.9 \\
%& 32.3 & 2.0 & 63.2 & 25.4 &  34.1 & 10.1 \\
%& 50.0 & & 17.4 & & 19.7 & \\
%\bottomrule &  &  &   
%\end{array}
%\end{equation*}%
%\caption{Results of the EM inversions using the SA and BFGS methods. The SNR is $\epsilon = 0.01.$}
%\label{tabella_risultati1}
%\end{table}

\begin{table}[tbp]
\begin{equation*}
\begin{array}{ccccccc}
\toprule
 & \multicolumn{2}{c}{\mathrm{\textbf{Mean} \; \textbf{relative} \; \textbf{error}} (SA)} & \multicolumn{2}{c}{\mathrm{\textbf{Mean} \; \textbf{relative} \; \textbf{error}} (BFGS)}  \\
 & \mathrm{Conductivity} & \mathrm{Thickness} & \mathrm{Conductivity} & \mathrm{Thickness} \\
 & (\%) & (\%) & (\%) & (\%) \\ 
\midrule
\multirow{3}*{Model 1} & 0.0 & 0.13 & 0.1 & 0.0 \\
& 11.5 & 1.86 & 5.4 & 5.4 \\
& 0.02 & & 0.1 & \\
\midrule
\multirow{3}*{Model 2} & 0.0 & 0.26 & 0.1 & 0.2 \\
& 3.38 & 4.51 & 2.0 & 1.1 \\
& 0.01 & & 0.1 & \\
\midrule
\multirow{3}*{Model 3} & 0.0 & 0.0 & 0.7 & 1.0 \\
& 8.71 & 7.1 & 53.0 & 23.8 \\
& 1.26 & & 1.8 & \\
\midrule
\multirow{3}*{Model 4} & 0.0 & 1.21 & 0.0 & 0.3 \\
& 6.66 & 15.9 & 1.7 & 6.2 \\
& 0.45 & & 0.5 & \\
\midrule
\multirow{1}*{Average error} & 2.66 & 3.87 & 5.4 & 4.75 \\
\bottomrule &  &  &   
\end{array}
\end{equation*}%
\caption{Results of the EM inversions using the SA and BFGS methods. The NSR is $\epsilon = 0.$}
\label{tabella_risultati1}
\end{table}

\begin{table}[tbp]
\begin{equation*}
\begin{array}{ccccccc}
\toprule
 & \multicolumn{2}{c}{\mathrm{\textbf{Mean} \; \textbf{relative} \; \textbf{error}} (SA)} & \multicolumn{2}{c}{\mathrm{\textbf{Mean} \; \textbf{relative} \; \textbf{error}} (BFGS)}  \\
 & \mathrm{Conductivity} & \mathrm{Thickness} & \mathrm{Conductivity} & \mathrm{Thickness} \\
 & (\%) & (\%) & (\%) & (\%) \\ 
\midrule
\multirow{3}*{Model 1} & 0.34 & 2.0 & 16.3 & 11.2 \\
& 24.7 & 33.4 & 24.2 & 5.5 \\
& 1.90 & & 9.2 & \\
\midrule
\multirow{3}*{Model 2} & 0.53 & 5.95 & 8.3 & 9.5 \\
& 22.1 & 19.7 & 7.6 & 14.8 \\
& 4.16 & & 2.6 & \\
\midrule
\multirow{3}*{Model 3} & 0.09 & 0.87 & 9.4 & 7.5 \\
& 19.0 & 19.2 & 15.2 & 21.0 \\
& 5.38 & & 3.4 & \\
\midrule
\multirow{3}*{Model 4} & 0.09 & 4.81 & 5.2 & 4.0 \\
& 23.1 & 20.2 & 1.8 & 7.0 \\
& 8.9 & & 6.2 & \\
\midrule
\multirow{1}*{Average error} & 9.19 & 13.3 & 9.12 & 10.0 \\
\bottomrule &  &  &   
\end{array}
\end{equation*}%
\caption{Results of the EM inversions using the SA and BFGS methods. The NSR is $\epsilon = 0.1 \%.$}
\label{tabella_risultati2}
\end{table}

\begin{table}[tbp]
\begin{equation*}
\begin{array}{ccccccc}
\toprule
 & \multicolumn{2}{c}{\mathrm{\textbf{Mean} \; \textbf{relative} \; \textbf{error}} (SA)} & \multicolumn{2}{c}{\mathrm{\textbf{Mean} \; \textbf{relative} \; \textbf{error}} (BFGS)}  \\
 & \mathrm{Conductivity} & \mathrm{Thickness} & \mathrm{Conductivity} & \mathrm{Thickness} \\
 & (\%) & (\%) & (\%) & (\%) \\  
\midrule
\multirow{3}*{Model 1} & 1.0 & 9.56 & 21.0 & 16.2 \\
& 32.1 & 34.9 & 35.2 & 4.5 \\
& 8.66 & & 12.6 & \\
\midrule
\multirow{3}*{Model 2} & 0.43 & 8.23 & 11.0 & 13.2 \\
& 27.4 & 30.1 & 21.5 & 20.7 \\
& 5.19 & & 4.5 & \\
\midrule
\multirow{3}*{Model 3} & 0.37 & 5.01 & 16.7 & 12.3 \\
& 24.3 & 31.7 &  29.1 & 19.2 \\
& 13.7 & & 4.8 & \\
\midrule
\multirow{3}*{Model 4} & 0.37 & 12.4 & 15.7 & 11.6 \\
& 32.2 & 24.4 &  13.7 & 8.5 \\
& 12.6 & & 18.4 & \\
\midrule
\multirow{1}*{Average error} & 13.2 & 19.5 & 17.02 & 13.28 \\
\bottomrule &  &  &   
\end{array}
\end{equation*}%
\caption{Results of the EM inversions using the SA and BFGS methods. The NSR is $\epsilon = 0.5 \%$}
\label{tabella_risultati3}
\end{table}

\section{Conclusions}

The numerical and analytical estimation of the electric and magnetic fields is at the basis of the modeling and inversion algorithms commonly used in induction EM methods. Here we present a novel methodology for the computation of the integrals involved in the solution of the Maxwell equations, based on the splitting of the reflection term. In this context, a classical quadrature rule on finite intervals can be applied. On the other hand, this approach also allows an analytical approximation of the integrals, that can be used to speed-up the solution of the inverse problems. \\
We have estimated the validity of our approach for the specific case of the DUALEM (DUAL-geometry Electro-Magnetic) system, but the physics is not restricted to these type of instruments.  
The specific examples consider the study of river-levees integrity, which is a very important environmental problem in Italy, due to the high hydrological risks. \\
A numerical algorithm based on the Gauss-Kronrod technique is used to compute the components of the EM field at low frequency in a stratified medium. The results are in good agreement with those of the commonly used digital filtering method. Moreover, the analytical approximation match well with the numerical solution both for resistive and conductive environments. \\ 
Two optimization algorithms are applied for the solution of the inverse problem, the line-search BFGS method, enhanced by the analytical approximation, and the global-search Simulated Annealing technique. The numerical experiments confirm the reliability of these techniques. Furthermore, while the two methods show similar performances in terms of solution accuracy, the former is more than 10 times faster than the latter. This asset of the enhanced BFGS method enables its application to the inversion of large datasets.  \\  
Although it is specific for the case of river levees, this analysis has a general validity, and allows to overcome the limits of common methods based on the modelling of apparent conductivity in the low induction number (LIN) approximation. 

\appendix

\section{Electrical conductivity}\label{AppA} 

The subsoil conductivities are computed with the complex refractive index model (CRIM).
The CRIM model for a shaly sandstone with negligible permittivity and partially saturated with gas, can be expressed as \cite{Schon, Carcione2014}
\begin{equation} \label{CRIM}
\sigma = [ (1- \phi) (1 - C) \sigma_q^{\gamma} + (1- \phi) C \sigma_c^{\gamma} + \phi S_w  \sigma_w ^{\gamma} +
\phi (1 - S_w) \sigma_a^\gamma ]^{1/\gamma}, \ \ \ \ \gamma = 1/2
\end{equation}
where $\sigma_q$, $\sigma_c$, $\sigma_w$ and $\sigma_a$ are the sand-grain (quartz), clay, water and air conductivities, $C$ is the clay content, $\phi$ is the porosity and $S_w$ is the water saturation. 
If $\gamma$ is a free parameter, the equation is termed Lichtnecker-Rother formula.
It is based on the ray approximation. The travel time in each medium is inversely proportional 
to the electromagnetic velocity, which in turn is inversely 
proportional to the square root of the complex dielectric constant. 
At low frequencies, displacement currents can be neglected and equation (\ref{CRIM}) is obtained. 
For zero clay content, and neglecting $\sigma_q$ and $\sigma_c$, equation (\ref{CRIM}) is exactly Archie's law \cite{Hoversten}.

\section*{Acknowledgements}

This work was partially supported by GNCS-INdAM, FRA-University of Trieste and CINECA under HPC-TRES program award number 2019-04. This research was also supported by the DILEMMA project (Imaging, modelling, monitoring and design of earthen levees), funded by the "Ministero dell’Ambiente e della Tutela del Territorio e del Mare" (MATTM). Eleonora Denich and Paolo Novati are members of the INdAM research group GNCS.


\begin{thebibliography}{99}

\bibitem[1]{Abramowitz} M. Abramowitz and I. A. Stegun (1970), \emph{Handbook of Mathematical Functions with Formulas, Graphs, and Mathematical Tables}, Dover Publications, Inc.,New York.

\bibitem[2]{anderson1} Anderson, W.L. (1982), \emph{Fast Hankel-transforms
using related and lagged convolutions}, ACM Transactions on Mathematical
Software, 8(4).

\bibitem[3]{anderson2} Anderson, W.L. (1979), \emph{Computer Program: Numerical integration of related 	Hankel transforms of orders 0 and 1 by adaptive digital filtering}, Geophysics, 44(7).

\bibitem[4]{Broyden} Broyden, C. G. (1970), \emph{The Convergence of a
Class of Double-Rank Minimization Algorithms}, Journal Inst. Math. Applic.

\bibitem[5]{Carcione2014} Carcione J. M. (2014), 
\emph{Wave fields in real media: Wave propagation in anisotropic, anelastic, porous and electromagnetic media}, Elsevier Science, Amsterdam, Handbook of Geophysical Exploration.

\bibitem[6]{Carcione2012} Carcione J. M., Gei D., Picotti S. and  Michelini A. (2012), 
\emph{Cross-hole electromagnetic and seismic modeling for CO$_2$ detection and monitoring in a saline aquifer}, Journal of Petroleum Science and Engineering.

\bibitem[7]{Erdelyi} Erdelyi, A., Ed. (1954), \emph{Tables of Integral
Transform}, McGraw-Hill Book Co.

\bibitem[8]{Fletcher} Fletcher, R. (1970), \emph{A New Approach to Variable
Metric Algorithms}, Computer Journal.

\bibitem[9]{Francese} Francese R., Morelli G., Santos F. M., Bondesan A., Giorgi M. and Tessarollo A. (2018), 
\emph{An integrated geophysical approach to scan river embankments}, Fast Times, 23(3), 86-96.

\bibitem[10]{Ghosh} Ghosh, D.P. (1971), \emph{The application of linear
filter theory to the direct interpretation of geoelectrical resistivity
sounding measuraments}, Geophysical Prospecting, 19(2), 192-217.

\bibitem[11]{Goffe} Goffe W. L., Ferrier G. D. and Rogers J. (1994), 
\emph{Global optimization of statistical functions with Simulated Annealing}, Journal of Econometrics
60(1). 

\bibitem[12]{Goldfarb} Goldfarb, D. (1970), \emph{A Family of Variable
Metric Updates Derived by Variational Means}, Mathematics of Computing, 24.

\bibitem[13]{GR} Gradshteyn, I. S. and Ryzhik, I. M. (1980), \emph{Tables of
Integrals, Series, and Products}, Academic Press, New York.

\bibitem[14]{Guptasarma} Guptasarma, D. (1982), \emph{Optimization of short
digital linear filters for increased accuracy}, Geophysical Prospecting 30,
501-514. 

\bibitem[15]{G-S} Guptasarma, D., Singh, B. (1997), \emph{New digital linear
filters for Hankel $J_0$ and $J_1$ transform}, Geophysical Prospecting 45.

\bibitem[16]{Hoversten} Hoversten G. M., Cassassuce F., Gasperikova E., Newman G. A., Chen J., Rubin Y., Hou Z. and Vasco D. (2006), 
\emph{Direct reservoir parameter estimation using joint inversion of marine seismic AVA and CSEM data.i}, Geophysics, 71(3), C1-C13.

\bibitem[17]{CR1} Ingeman-Nielsen, T. and Baumgartner, F. (2006), \emph{%
CR1Dmod: A Matlab program to model 1D complex resistivity effects in
electrical and electromagnetic surveys}, Computer \& Geosciences, 32,
1411-1419. 

\bibitem[18]{JS} Johansen, H. K., and Sorensen, K. (1979), \emph{Fast Hankel
transforms}, Geophysical Prospecting, 27(4), 876-901. 

\bibitem[19]{Kirkpatrick} Kirkpatrick S., Gelatt C. D. and Vecchi M. P. (1983), 
\emph{Optimization by Simulated Annealing}, Science
4598(220). 

\bibitem[20]{Wiener Hopf} Koefoed, O., and Dirks, F.J.H. (1979), \emph{%
Determination of resistivity sounding filters by the Wiener-Hopf
least-squares method}, Geophysical Prospecting, 27(1), 245-250.

\bibitem[21]{KGP} Koefoed O., Ghosh, D. P., and Polman, G. J. (1972), \emph{%
Computation of type curves for electromagnetic depth sounding with a
horizontal transmitting coil by means of a digital linear filter},
Geophysical Prospecting, 20(2), 406-420.

\bibitem[22]{Kong} Kong, F.N. (2007), \emph{Hankel transform filters for
dipole antenna radiation in a conductive medium}, Geophysical Prospecting.

\bibitem[23]{Pawar} Pawar S. D., Murugavel P. and Lal D. M. (2009), 
\emph{Effect of relative humidity and sea level pressure on electrical conductivity of air over Indian Ocean}, Journal of Geophysical Research, 114(D2).

\bibitem[24]{Schon} Sc{h\"o}n J.H. (1996), 
\emph{Physical properties of rocks}, Pergamon Press, Handbook of Geophysical Exploration.

\bibitem[25]{quadgk} Shampine, L.F. (2008), \emph{Vectorized Adaptive Quadrature in MATLAB$^{®}$}, Journal of Computational and Applied Mathematics, 211, 131-140. 

\bibitem[26]{Shanno} Shanno, D. F. (1970), \emph{Conditioning of
Quasi-Newton Methods for Function Minimization}, Mathematics of Computing.

\bibitem[27]{Emdpler} Singh, N.P., Mogi, T. (2010), \emph{EMDPLER: A F77
program for modelling the EM response of dipolar sources over the
non-magnetic layer earth models}, Computer \& Geosciences, 36, 430-440.

\bibitem[28]{Palacky} Ward S. H. and Hohmann G. W. (1987), 
\emph{Resistivity characteristics of geologic targets}, In: Nabighian, M.N., (ed.),
  Electromagnetic Methods in Applied Geophysics, Tulsa, Oklahoma, 53-129.

\bibitem[29]{WH} Ward, S.H. and Hohmann, G.W. (1988), \emph{Electromagnetic
theory for geophysical applications}, In: Nabighian, M.N., (ed.),
Electromagnetic Methods in Applied Geophysics, Tulsa, Oklahoma, pp 131-311.

\bibitem[30]{Bessel} Watson, G. N. (1966), \emph{A treatise on the theory of Bessel functions}, Cambridge University press, 2nd. edition reprinted.

\bibitem[31]{Werthmuller} Werthmuller, D., Key, K., and Slob, E.C. (2018), 
\emph{A tool for designing digital filters for the Hankel and Fourier
transforms in potential, diffusive, and wavefield modeling}, Geophysics.

\bibitem[32]{Zohdy} Zohdy A. A. R. and Jackson D. B. (1969), 
\emph{Application of deep electrical soundings for groundwater exploration in Hawaii}, Geophysics, 34(4), 584-600.





\end{thebibliography}
\end{document}